\documentclass{amsart}

\usepackage{amsthm,amsfonts,amsmath,amssymb,latexsym,epsfig}

\usepackage{upref,amssymb,eucal,ae}

\newtheorem{theorem}{Theorem}
\newtheorem{lemma}[theorem]{Lemma}
\newtheorem{proposition}[theorem]{Proposition}
\newtheorem{corollary}[theorem]{Corollary}
\newenvironment{example}{\medskip \refstepcounter{theorem}
\noindent  {\bf Example \thetheorem}.\rm}{\,}
\newenvironment{remark}{\medskip \refstepcounter{theorem}
\noindent  {\bf Remark \thetheorem}.\rm}{\,}

\renewcommand{\thetheorem}{\thesection.\arabic{theorem}}

\def\b#1{\overline{#1}}
\def\ddb{\partial  \overline{\partial}}

\def\bn#1#2{\left( \begin{array}{c} #1 \\ #2 \end{array} \right)}

\def\tg{{\tilde g}}

\def\<{\langle}

\def\n{\nabla}
\def\>{\rangle}

\def\a{\alpha}

\def\o{\omega}

\def\b{\beta}
\def\tm{\tilde{M}}
\def\tn{\tilde{n}}
\def\tg{\tilde{g}}
\def\mb#1{{\mathbb #1}}
\def\mc#1{{\mathcal #1}}
\def\mf#1{{\mathfrak #1}}
\begin{document}

\title[The second fundamental form of isometric immersions]{The $L^2$-norm 
of the second fundamental form of isometric immersions into a Riemannian manifold}
\author{Santiago R. Simanca}
\email{srsimanca@gmail.com}

\begin{abstract}
We consider critical points of the global squared $L^2$-norms of the second 
fundamental form and the mean curvature vector of isometric immersions into 
a fixed background Riemannian manifold under deformations of the immersion. 
We use the critical points of the former functional to define 
canonical representatives of a given integer homology class of the background 
manifold. We study the fibration ${\mathbb S}^3 \hookrightarrow Sp(2)\stackrel 
{\pi_{\circ}}{\rightarrow} {\mathbb S}^7$ from this point of view, showing 
that the fibers are the canonical generators of the $3$-integer homology of 
$Sp(2)$ when this Lie group is endowed with a suitable family of left 
invariant metrics.  Complex subvarieties in the standard $\mb{P}^n(\mb{C})$ 
are critical points of each of the functionals, and are canonical
representatives of their homology classes. We use 
this result to provide a proof of Kronheimer-Wrowka's theorem on the smallest 
genus representatives of the homology class of a curve of degree $d$ in 
${\mathbb C}{\mathbb P}^2$, and analyze also the canonical representability of 
certain homology classes in the product of standard $2$-spheres. Finally, 
we provide 
examples of background manifolds admitting isotopically equivalent critical 
points in codimension one for the difference of the two functionals 
mentioned, of different critical values, which are Riemannian analogs of 
alternatives to compactification theories that has been offered recently.  
\end{abstract}

\maketitle

\section{Introduction}
Consider a fixed background Riemannian manifold $\tm$ of a given dimension, 
and let $M$ be an immersed closed submanifold.
With the metric induced from that of $\tm$, $M$ itself becomes a Riemannian 
manifold. Of all such isometric immersions of $M$ into $\tm$, we would like 
to single out the one that is {\it optimal}. This is the main theme of our
work here.

Before going any further, let us relate the setting
just described to the language of contemporary mathematical physics \cite{pe}. 
We liken $\tm$ to the full world, the {\it bulk} in the
terminology of physicists. The background metric it carries makes it possible
to define the length of curves in the ambient space, and consequently, we
can talk of the {\it least action principle} for the motion of particles in
this world. In this bulk, the isometrically embedded submanifold $M$ is taken
to represent 
a {\it brane}. Branes seem 
to come in two flavours, those bounding portions of the ambient space, or those
serving as the boundary locations for open strings. Our setting differs from
that in physics in that the ambient space metric is Riemannian, and not
Lorentzian. But we dispense with that defect in the analogy, and move along
anyway. The question of interest
to us can be phrased by asking for the optimal way of isometrically
immersing a brane into the bulk.

An extremely worthy effort in this direction can be made if we think of the
optimal immersion of a brane as one that extremizes the volume the brane
occupies among deformations of the immersion. This of course leads to the
study of minimal submanifolds, which at least in the case of surfaces goes
back to Lagrange, or even a bit earlier, to Euler, who studied the catenoid, 
in the mid 1700s, and which constitutes a subject that have 
remained an exciting area of research to this date. 

Here, we turn to a different criteria in order
to measure optimality, again motivated by the development in theoretical
physics. Indeed, it seems that the {\it Calabi-Yau} manifolds, a particular
type of K\"ahler manifold, are to be the preferred choice for compactifying or
curling-up the extra dimensions we do not seem to be able to observe in our
bulk world. And though all complex submanifolds of a K\"ahler manifold are 
minimal, there are too many minimal submanifolds that are not complex. The
said volume criteria does not give a special role to K\"ahler branes, and
in that sense, there might be others that are more appropriate in order 
to measure the optimality of an immersion, perhaps singling out 
K\"ahler branes over other minimal submanifolds. 

Since our submanifold $M$ is isometrically immersed into the ambient space
$\tm$, we can measure the length of paths whose end points are in $M$ once we 
know the immersion.  Thus, the issue is strictly about how to place the 
submanifold $M$ inside the ambient space $\tm$. Since the metric in
$\tm$ defines canonically a notion of derivative everywhere,
the Levi-Civita connection of its metric, we may naturally conceive of an energy 
associated to an immersion that only depends upon its first order derivatives,
and so as not to alter the submanifold itself, make this energy depend 
upon the normal component of these first derivatives. Thus, we look at the
second fundamental form of the immersion, $\alpha$, and declare that the 
immersion's energy be measured by its global squared $L^2$-norm, which we shall
call $\Pi (M)$. This energy for an immersion resembles 
the stored-energy function for neo-Hookian materials, and by way of analogy, we
may think of the energy of an immersed brane as being determined by the
accumulated effect of a tensor of {\it springs} that act in the 
normal directions. In accordance with Hamilton's principle, an optimal 
isometric immersion of the submanifold will be one that is a stationary point
for $\Pi (M)$ under deformations of the isometric immersion. Since our 
immersion does not have {\it kinetic energy}, this
criterion for the optimality of the immersion corresponds exactly to that used
to place a static elastodynamic incompressible neo-Hookian material in space,
the incompressibility taken here as the analogue of an immersion that is
isometric. 
Notice that in this sense, the functional $\Pi(M)$ is a much 
{\it higher order energy} than that given by the volume of the immersed 
manifold.

We develop the Euler-Lagrange equations for $\Pi(M)$ 
under deformations of the isometric immersion of $M$ into $\tm$.
As equations on the coefficients of the second fundamental form of an 
stationary immersion, they constitute a {\it local} nonlinear system of order two,
with cubic nonlinearities in the principal curvatures along normal directions.
(The local nature of the equations shows at least one defect with the analogy 
with the motion of 
incompressible elastic bodies above, which are ruled by pseudodifferential 
equations instead.)    
In the case of codimension one, we conclude that all austere manifolds 
are critical points. Though the set of these is much richer in general, 
it is likely that the critical points of $\Pi(M)$ characterize austerity in 
dimensions less or equal than $4$. We also prove that all complex 
submanifolds of complex projective space with the Fubini-Study metric 
are stationary points, in a sense
vindicating the desired preference that we wanted to give to K\"ahler branes. 

For a given homology class of dimension $k$, the embedded minimizers of 
$\Pi(M)$ within the class are to be considered its canonical 
representatives.  In contrast with the case where $k=1$, these minimizers do 
not have to be totally geodesics, or even exist for that matter. We prove their
existence in the case where $\tm=Sp(2)$ endowed with a suitable left-invariant
metric, and the homology class is the generator of $H_3(Sp(2);\mb{Z})$. This 
same background manifold can be used to show the nonexistence
of embedded minimizers of $\Pi(M)$ representing the generator of 
$H_7(Sp(2);\mb{Z})$, for the simple reason that this class does not
admit embedded representatives.

The cubic nature of the nonlinearity in the Euler-Lagrange equation for 
$\Pi(M)$ seems to single out K\"ahler branes of complex dimension less 
or equal than $3$. Exploting the fact that complex submanifolds are minimal 
(see, for instance, \cite{bl,si}), we refine our criterion for optimality of
the immersion, and take it one step further.
The trace of the second fundamental form of an immersion 
is the mean curvature vector, $H$, whose squared $L^2$-norm we shall 
call $\Psi(M)$,
the Willmore functional \cite{wi}. 
We develop also its Euler-Lagrange equation under deformations of the isometric
immersion, and having done so, consider the ``energy'' functional 
$\Pi(M)-\Psi(M)$ over the space of deformations of the isometric immersion 
restricted to represent a fixed integral homology class. In a sense, 
stationary points of this functional
are given by isometric immersions for which the trace-free part of the
second fundamental form is as small as possible, paralleling the defining
property of Einstein metrics, whose Ricci tensors are traceless.
Thus, we look at these stationary points as some sort of Einstein 
isometric immersions of the brane into the bulk. If $M$ is a canonical
isometric immersion of the smallest volume among critical points of
$\Pi(M)$ of a fixed critical level, when $M$ itself is minimal, it is then also 
a critical point of $\Pi(M)-\Psi(M)$, and in fact, that stationary metric
on $M$ is a critical point of the total scalar curvature functional in the
space of Riemannian metrics on $M$ that can be realized by isometric 
embeddings into $\tm$.    

Complex submanifolds of complex projective space are examples of such.
For when the ambient space $\tm$ is K\"ahler, any complex
submanifold is an absolute minimum for $\Psi(M)$, and so the said
complex submanifolds are critical points of $\Pi(M)-\Psi(M)$ under
deformations. Thus, the new functional distinguishes complex submanifolds 
at least when the bulk background is equal to complex projective space endowed 
with the Fubini-Study metric. We use this property, and a description of the
singularities of minimizing sequences of representatives of a given 
integral homology class, to prove that among all embedded surfaces in 
complex projective plane that represent the homology class
of a complex curve of fixed degree, the complex ones are the ones with the
smallest genus. This is a theorem of Kronkheimer and Mwroka \cite{krmr},
after a conjecture commonly attributed to R. Thom.

The energies $\Pi(M)$, $\Pi(M)-\Psi(M)$, and $\Psi(M)$   
differ from each other, 
and this is manifested well in the properties of their stationary points. The
differences remains in place when our functionals are compared to others 
that have been termed Willmore in the literature
also \cite{glw,hl}. For in the case 
of codimension one, for instance, Example \ref{ex7}
below exhibits a family of 
Willmore hypersurfaces out of which $\Pi(M)$ singles out as stationary 
points only those that are Clifford minimal tori. 

The choice of $\Pi(M)-\Psi(M)$ as an energy functional for the immersion is
very natural, and truly derives from equation (\ref{sc}) below that expresses 
the intrinsic scalar curvature of $M$ in terms of the scalar curvature of 
$\tm$ and additional extrinsic data, an identity that in turn follows from 
Gauss' equation. 
We obtain $\Pi(M)-\Psi(M)$ by integrating over the immersed manifold 
the summands in this expression that involve first order covariant derivatives 
of the connection of the ambient space metric.  Thus, the said functional is 
the averaged global extrinsic first order contribution
to the total scalar curvature of the metric on $M$ induced by
that of $\tm$. The remaining unaccounted contribution to the intrinsic total 
scalar curvature is obtained by averaging the sectional curvatures of the 
ambient space metric over the submanifold, and in that sense, it is not a 
contribution that we are at liberty to move by changing the isometric immersion.  
Its minimal critical points are thus critical points of the total scalar
curvature functional over the space of metrics that can be realized by 
isometric immersions into $\tm$s, and many of these, are Einstein
metrics on $M$, cf. \cite{rss}.  

\subsection{Structure of the article}
In \S\ref{s2} we introduce the basic notation, and conventions. As we use them 
later on, we derive expressions for the Ricci tensor and scalar curvature of
an immersed manifold in terms on intrinsic and extrinsic data.
In \S\ref{s3} we derive the Euler-Lagrange equation for $\Pi(M)$ and $\Psi(M)$ 
under deformations of the immersion. We write these as a system of $q$ 
equations, $q$ the codimension of $M$ in $\tm$, and discuss the case
where $M$ is a hypersurface and the background metric is Einstein. We 
consider also the functional $\Pi(M)$ when its domain is restricted to the set of 
embedded representatives of a fixed integer homology class $D$ of $\tm$, and 
define a canonical representative of $D$ to be a minimum $M$ of smallest
volume.  The theory is illustrated using the simple Lie group
$Sp(2)$ and the generator $D$ of the third homology group $H_3(Sp(2);
{\mathbb Z})$.
We prove that if $Sp(2)$ is endowed with suitable left-invariant metrics, the
canonical representative of $D$ exists and, up to isometric diffeomorphisms,
is unique. 
The special nature of complex submanifolds of complex projective space
is presented in \S\ref{s4}, where we briefly recall also some relevant facts 
about the Fubini-Study metric in general, and curves in complex projective 
plane, both to be used in later portions of our work. We prove that if $(M,g)$ is
a complex manifold isometrically immersed into $\mb{P}^n(\mb{C})$ with its Fubini-Study metric, then $\| \alpha \|^2$ is an intrinsic quantity, and develop various
other properties of complex subvarieties in relation to the functional 
$\Pi(M)$.  We continue the analysis of these properties  
and its relation to canonical representatives of homology classes in  
\S\ref{s5}. In combination with the fact that complex curves minimize 
$\Psi(M)$, we derive in there an alternative proof of
the aforementioned theorem of Kronkheimer and Mwroka,  
avoiding the use of the Seiberg-Witten theory.
In this section, we discuss also the canonical representability of the class
$(1,1)$ in the 2-dimensional homology of the product of two spheres.
Finally, in \S\ref{s6}, we exhibit examples of critical hypersurfaces of
$\Pi(M)-\Psi(M)$ of different critical values that are 
isotopically equivalent to each other in a) ${\mathbb S}^{n+1}$,
$n\geq 8$, provided with the standard metric; b) the trivial bundle 
${\mathbb S}^2\times \mb{R}^2$ provided with the Schwarzchild
metric; and c) the bundle
$T{\mathbb S}^2$ provided with the Eguchi-Hanson metric. In the last
two cases, the second of the exhibited 
critical hypersurfaces occurs asymptotically
at $\infty$.

\section{Preliminaries, notations and conventions}
\label{s2}
\setcounter{theorem}{0}
Consider a closed Riemannian manifold $(\tm, \tg)$, and let $M$ be a 
submanifold of $\tm$. With the induced metric, $M$ itself becomes 
a Riemannian manifold. We denote this metric on $M$ by $g$, and proceed to
study certain relations between the geometric quantities of $g$ and $\tg$ on
points of $M\subset \tm$. We shall denote by $n$ and $\tn$ the dimensions
of $M$ and $\tm$, respectively.

The {\it second fundamental form} $\a$ of $M$ is defined in terms of the
decomposition of the Levi-Civita connection $\n^{\tg}$ into a tangential and 
normal component \cite{md}. Indeed, under the inclusion map, the pull-back 
bundle $i^{*}T\tm$ can be decomposed as a Whitney sum
$$i^{*}T\tm = TM \oplus \nu(M)\, ,$$
where $\nu(M)$ is the normal bundle to $M$. Given vector fields $X$ and $Y$
tangent to $M$, we have that 
$$
\a(X,Y)=\pi_{\nu(M)}(\n^{\tg}_X Y )\, ,
$$
where $\pi_{\nu(M)}$ is the orthogonal projection onto the normal bundle. 

By uniqueness of the Levi-Civita connection, the orthogonal projection of the
ambient space connection onto $TM$ defines $\n^{g}_X Y$. We have Gauss' 
identity
\begin{equation}
\n^{\tg}_X Y = \n^{g}_X Y + \a(X,Y) \, . \label{ga}
\end{equation}

If we consider a section $N$ of the normal bundle $\nu(M)$, the shape operator 
is defined by 
$$
A_N X= - \pi_{TM}(\n^{\tg}_X N )\, ,
$$
where in the right side above, 
$N$ stands for an extension of the original section to
a neighborhood of $M$. 
 If $\n^{\nu }$ is the connection on $\nu(M)$ induced by $\n^{\tg}$, we have
Wiengarten's identity
\begin{equation}
\n^{\tg}_X N = -A_N X + \n^{\nu}_X N \, .\label{we}
\end{equation}

The identities (\ref{ga}) and (\ref{we}) allows us to express the 
commutator $[\n^{\tg}_X,\n^{\tg}_Y]$ in terms of the commutator 
$[\n^{g}_X,\n^{g}_Y]$,
the second fundamental form and the shape operator. We may do so also 
for  $\n^{\tg}_{[X,Y]}$. As a consequence, if $Z$, $W$ are two additional 
vector fields tangent to $M$, we obtain Gauss' equation, a fundamental 
relation 
between the tangential component of the extrinsic and intrinsic 
curvature tensors on $M$:
\begin{equation}
g(R^g(X,Y)Z,W)= \tg(R^{\tg}(X,Y)Z,W)+\tg(\a(X,W),\a(Y,Z))-\tg(\a(X,Z),\a(Y,W))
\, .\label{cu}
\end{equation}
Here, $R^g$ stands for the Riemann curvature tensor of the corresponding 
metric $g$.

The Ricci tensor $r_g(X,Y)$ of a Riemannian metric $g$ with
curvature tensor $R^g$ is defined as the trace of the linear map 
$L\rightarrow R^g(L,X)Y$, while the scalar curvature $s_g$ is defined as the
metric trace of $r_g$.
Let $\{e_1, \ldots , e_{\tn}\}$ be an orthonormal frame for $\tg$ in a 
neighborhood of a point in $M$ such that $\{e_1, \ldots , e_{n}\}$ 
constitutes an orthonormal frame for $g$ on points of $M$. Letting $H$ 
be the mean curvature
vector, the trace of the second fundamental form, by (\ref{cu}) we obtain
that 
$$\begin{array}{rcl}
r_g(X,Y) \! \! \! & = & \! \! \! \sum_{i=1}^n \tg(R^{\tg}(e_i,X)Y,e_i)+
\tg(H,\a(X,Y))-\sum_{i=1}^n \tg(\a(e_i,X),\a(e_i,Y)) \\ & = & 
r_{\tg}(X,Y) - \sum_{i=n+1}^{\tn} \tg(R^{\tg}(e_i,X)Y,e_i) +\tg(H,\a(X,Y))-
 \\ & & \mbox{} \hspace{2in} \sum_{i=1}^n \tg(\a(e_i,X),\a(e_i,Y)) \, ,
\end{array}$$
while
\begin{equation}
\begin{array}{rcl}
s_g & = & s_{\tg}-2 \sum_{i=1}^n \sum_{j=n+1}^{\tn}K_{\tg}(e_i,e_j)
-\sum_{i,j>n}K_{\tg}(e_i,e_j) +\tg(H,H)- \tg(\a,\a) \vspace{1mm} \\
& = & \sum_{i,j\leq n}K_{\tg}(e_i,e_j)
+\tg(H,H)- \tg(\a,\a) \, .
\label{sc}
\end{array}
\end{equation}
Here $K_{\tg}(e_i,e_j)$ is the $\tg$-curvature of the section spanned by
the orthonormal vectors $e_i$ and $e_j$, and $\tg(H,H)$ and $\tg(\a,\a)$ are
the squared-norms of the mean curvature vector $H$ and the form $\a$,
respectively.

In the same manner as Gauss' equation (\ref{cu}) describes the tangential 
component of 
$R^{\tg}(X,Y)Z$ in terms of the intrinsic curvature tensor of $M$ and the
second fundamental form, we may find an analogous type of relation 
for the normal component $(R^{\tg}(X,Y)Z)^{\nu}$ of this same vector. In so
doing, we obtain the following identity, known as Codazzi's equation: 
\begin{equation}
\begin{array}{rcl}
\n_{X}^{\tg}\a(Y,Z) & = & (R^{\tg}(X,Y)Z)^{\nu}+\n^{\tg}_{Y}\a(X,Z)+
A_{\a(X,Z)}Y- A_{\a(Y,Z)}X \\ &  & \mbox{}\hspace{0.2in}
\a(Y,\n_X^g Z)- \a(X,\n_Y^g Z)
+ \a([X,Y],Z) \, .
\end{array}
\label{co}
\end{equation}
If $\{ e_i\}_{i=1}^n$ is an orthonormal frame of $TM$ as above that is
geodesic at a point $p$, we derive the following identity
$$
\n_{e_j}^{\tg}\a(e_i,e_j)=(R^{\tg}(e_j,e_i)e_j)^{\nu}+\n^{\tg}_{e_i}H+
A_{H}e_i-A_{\a(e_i,e_j)}e_j \, , 
$$
which will be used several times below. In writing this expression, we tacitly
use the standard summation convention over repeated indexes, a practive to be 
repeated throughout the article. 

\section{Variational formulas} \label{s3}
\setcounter{theorem}{0}

We now develop the Euler-Lagrange equations of the squared $L^2$-norm of 
the second fundamental form, mean curvature vector, and total
extrinsic scalar curvature of a manifold $M$ 
isometrically immersed into $(\tm,\tg)$ under deformations of the immersion.
Since the densities of the functionals in question are {\it local} 
operators, in deriving the said equations we may
assume that the immersion is in fact a compact embedding. In what follows, 
the ambient space metric $\tg$ will be denoted by 
$\< \, \cdot \, , \, \cdot \, \>$, for convenience.

\subsection{The second fundamental form}
We consider a closed manifold $M$ of dimension $n$, and let 
$f:M \rightarrow \tm$ be an isometric immersion into the
$\tn$-dimensional Riemanniann manifold $(\tm,\tg)$. For the reasons
indicated above, we shall assume that $f(M)=M\hookrightarrow \tm$ is a compact 
embedding. Let $\a$ be the second fundamental form induced by $f$,
and $\| \a\|^2=\< \a, \a\>$ be its squared pointwise $L^2$-norm, 
expression which of course involves the metric on $M$. 
We let $\Pi(M)$ be the integral of $\| \a\|^2$. That is to
say, if $d\mu=d\mu_M$ denotes the volume measure on $f(M)$, then  
\begin{equation}
\Pi(M)= \int_{M} \| \a \|^2 \, d\mu \, . \label{sf}
\end{equation}

Let $f: (-a,a)\times M \rightarrow \tm$ be a one parameter family of 
deformations of $M$. We set $M_t = f(t,M)$ for $t\in (-a,a)$, and have 
$M_0=M$. We would like to compute the $t$-derivative of 
$\Pi(M_t)= {\displaystyle \int_{M_t} \| \a_t\|^2\, d\mu_t} $ at $t=0$. 

Given a point $p\in M$, we let $\{x^1, \ldots, x^n,t\}$ be a coordinate
system of $M\times (-a,a)$, valid in some neighborhood of $(p,0)$,
such that $\{x^1, \ldots, x^n \}$ are normal coordinates of $M$ 
at $p$. 
The induced metric $g_t$ on $M_t$ has components 
$g_{ij}=\< e_i, e_j\>$, where 
$e_i=df(\partial_{x^i})$, $i=1, \ldots, n$, and by assumption, 
$g_{ij}(p,0)=\delta_{ij}$. We set $T=df(\partial_{t})$, the variational 
vector field of the deformation. 

We decompose $T$ into its tangential and normal component, $T=T^t+T^n$.
Then the variation of the measure is given by
\begin{equation}
\frac{d\mu_t}{dt}(p,0)=\left( {\rm div}(T^t) - \< T^n,H\>\right)
d\mu_0(p,0) \, , \label{vf}
\end{equation}
a well-known expression (see, for instance, \cite{li}) 
already used extensively in the theory of minimal 
submanifolds\footnotemark.
\footnotetext{Formula (\ref{vf}) implies that
$M$ is a critical point of the volume functional 
under deformations of the isometric immersion if, and only if, $H=0$. 
Manifolds with this property are called {\it minimal}, as in our 
Introduction here, \S 1.}    
Since the gradient of a function defined on $M$ is a tangential vector 
field while the mean curvature is normal, we have that 
\begin{equation}
\dot{\Pi}(M_t)\mid_{t=0}=\frac{d}{dt}\int_{M_t}\|\a_t \|^2 \, d\mu_t 
\mid_{t=0}=\int_{M} \left( \frac{d}{dt} \| \a_t \|^2 \mid_{t=0}
- \< \n \| \a \|^2 +\| \a \|^2 H,T\> \right) d\mu \, .\label{var}
\end{equation}

When computing the variation of the pointwise $L^2$-norm of $\a$, it is
necessary to use an orthonormal frame along points $(p,t)$ of an integral 
curve of $T$ for small $t$s. The coordinates $\{ x^1, \ldots, x^n\}$ above 
were chosen so that the coordinate vector fields 
$\{ \partial_{x^j}\}$ are orthonormal and have vanishing $g$-covariant 
derivatives at $p$. These vector fields commute with $T$. 
We now take an orthonormal tangent frame $\{ e_1, \ldots, e_n\}$ along a 
sufficiently small neighborhood of the integral curve of $T$ through $p$ that 
extends the $\partial_{x^j}\! |_{p}$s. Thus, 
$t$ parametrizes the integral curve through 
$p$, and for each sufficiently small $t$, we have 
that $\{ e_1, \ldots, e_n\}$ is a normal frame of $M_t$ at $(p,t)$. We still
have that $[T,e_i]=0$ and $\n^{g}_{e_j}e_i=0$ at $p=(p,0)$, respectively.
With the Einstein summation convention in place, $\| \a\|^2$ can
be conveniently calculated by 
$\| \a \|^2=\< \a(e_i,e_j),\a(e_i,e_j)\>$, as already tacitly 
done in (\ref{sc}).

By differentiation along the integral curves of $T$, we have that 
\begin{equation}
\begin{array}{rcl}
{\displaystyle \frac{d}{dt}\< \a (e_i, e_j),\a(e_i,e_j)\>\mid_{t=0}} & = & 
2\< \n_T^{\tg}\a(e_i,e_j),\a(e_i,e_j)\>+2\dot{g}^{il}\< \a (e_i, e_j),
\a(e_l,e_j)\> \\ & = & 
2\< \n_T^{\tg} ( \n^{\tg}_{e_i}{e_j} -\n^{g}_{e_i}{e_j}) ,\a(e_i,e_j)\>
+2\dot{g}^{il}\< \a (e_i, e_j),
\a(e_l,e_j)\> \, .
\end{array} \label{pwd}
\end{equation}

Since $\n^{g_t}_{e_i}e_j$ is tangent to $M_t$, and vanishes at $p$, while
$\a(e_i,e_j)$ is a normal vector, we have that
$$\< \n_T^{\tg}\n^{g}_{e_i}{e_j} ,
\a(e_i,e_j)\> = T\< \n^{g}_{e_i}{e_j} , \a(e_i,e_j)\>
-\<\n^{g}_{e_i}{e_j} ,\n^{\tg}_T\a(e_i,e_j)\>=0\, .$$
On the other hand, at $p$ we have that
$$\begin{array}{rcl}
\< \n_T^{\tg} \n^{\tg}_{e_i}{e_j}, \a(e_i,e_j)\> & = &   
\< \n^{\tg}_{e_i}\n^{\tg}_{T}e_j +\n^{\tg}_{[T,e_i]}e_j+R^{\tg}(T,e_i)e_j,
 \a(e_i,e_j)\> \vspace{1mm} \\ & = & 
\< \n^{\tg}_{e_i}\n^{\tg}_{T}e_j +R^{\tg}(T,e_i)e_j,
 \a(e_i,e_j)\>
\end{array}$$
because $\n^{\tg}_X Y$ is tensorial in $X$, and $[T,e_i]_p=0$. 

Still computing at $p$, we find that
$$\begin{array}{rcl}
\< \n_{e_i}^{\tg} \n^{\tg}_{T}{e_j}, \a(e_i,e_j)\> & = &   
\< \n^{\tg}_{e_i}\n^{\tg}_{e_j}T +\n^{\tg}_{e_i}[T,e_j],
 \a(e_i,e_j)\>  \\ & = & \< \n^{\tg}_{e_i} \n^{\tg}_{e_j} T, 
\a(e_i,e_j)\> + e_i\< [T,e_j],\a(e_i,e_j)\> \\ & = &  
\< \n^{\tg}_{e_i} \n^{\tg}_{e_j} T, 
\a(e_i,e_j)\> \, ,
\end{array}
$$
with the second and third equality being true because we have that $[T,e_j]_p=0$, and with
the latter of these making use of the tangentiality of $e_i$ also. Further, 
$$
\< \n^{\tg}_{e_i}\n^{\tg}_{e_j}T, \a(e_i,e_j)\> \! = \! 
e_i e_j \< T, \a(e_i,e_j) \>-e_i\<T,\n^{\tg}_{e_j}\a(e_i,e_j)\>-e_j\<T,\n^{\tg}_{e_i}
\a(e_i,e_j)\> +
\< T, \n_{e_j}^{\tg} \n_{e_i}^{\tg}\a(e_i,e_j)\> \, ,
$$
expression whose right hand side is equal to its very last term 
modulo a differential, and so by Stokes' theorem, its integral over $M$ 
coincides with the integral of the said term.

We use the symmetries of the curvature tensor to obtain the identity 
$$
\< R^{\tg}(T,e_i)e_j,\a(e_i,e_j)\> =\<R^{\tg}(\a (e_i,e_j),e_j)e_i, T\>\, .
$$

Finally, using the orthonormality of the frame, we conclude that
$$
\dot{g}^{il}\< \a (e_i, e_j),
\a(e_l,e_j)\>=2\<\< \a (e_i, e_j),
\a(e_l,e_j)\> \a(e_i,e_l),T \> -(e_i\< T,e_l\> +e_l\<T,e_i\>)
\< \a (e_i, e_j),\a(e_l,e_j)\> \, .
$$

\begin{theorem} \label{ele}
Let $f: (-a,a)\times M \rightarrow \tm$ be a deformation of 
an isometrically immersed submanifold $f: M\rightarrow f(M)\hookrightarrow 
\tm$ into the Riemannian manifold $(\tm, \tg)$. We set $M_t=f(t,M)$, and have
$M_0=M$. We let $\{ e_1, \ldots , e_n\}$ be an orthonormal frame of $M_t$ 
for all $t\in (-a,a)$. Then the infinitesimal variation of {\rm (\ref{sf})} is 
given by
$$
\begin{array}{rcl}
{\displaystyle \frac{d \Pi(M_t)}{dt}\mid_{t=0}} & = & 
{\displaystyle \int_{M} \< 2\n_{e_j}^{\tg}\n_{e_i}^{\tg}\a(e_i,e_j)+
2R^{\tg}(\a (e_i,e_j),e_j)e_i
 -\| \a \|^2 H,T\> d\mu } \vspace{1mm} \\
& & + {\displaystyle \int_{M} \< 4\< \a(e_i,e_j),\a(e_l,e_j)\>
\a(e_l,e_i)-\n \| \a \|^2 , T\>d\mu 
} \vspace{1mm} \\
& & -2{\displaystyle \int_{M} (e_i\< T,e_l\> +e_l\<T,e_i\>)
\< \a (e_i, e_j),\a(e_l,e_j)\>d\mu }
\, .
\end{array}
$$
The isometrically immersed submanifold 
$f: M\rightarrow f(M)\hookrightarrow \tm$ satisfies the Euler-Lagrange 
equation of {\rm (\ref{sf})} if, and only if, the normal
component of the vector 
\begin{equation}
2\n_{e_j}^{\tg}\n_{e_i}^{\tg}\a(e_i,e_j) +2R^{\tg}(\a (e_i,e_j),e_j)e_i
 - \| \a \|^2 H + 4\< \a(e_i,e_j),\a(e_l,e_j)\>
\a(e_l,e_i)
\label{el}
\end{equation}
is identically zero.
\end{theorem}

{\it Proof}. The calculations preceding the statement may be 
used to derive a convenient expression for (\ref{pwd}). When this expression is
inserted into (\ref{var}), the variational formula results by applying 
Stokes' theorem.

We now compute this variational expression after decomposing the vector 
field $T$ into its tangential and normal components. Since a tangential 
deformation merely yields a reparametrization of $M$, and since 
the vector $\n \| \a\|^2$ is tangent to this manifold, the 
Euler-Lagrange equation of (\ref{sf}) is equivalent to the vanishing
of the normal component of (\ref{el}), as 
stated. \qed

We shall say that $(M,g)$ is canonically placed in $(\tm,\tg)$ if
there exists an isometric embedding $i: (M,g) \hookrightarrow (\tm,\tg)$ 
that is a critical point of (\ref{sf}) under deformations of the
embedding. There is no rigidity in how one of these manifolds can be placed 
inside $\tm$. We exhibit in \S \ref{s37} a family of isometric embeddings
$i_t: (M,g)\hookrightarrow (\tm,\tg)$
parametrized by a $7$-sphere, such that $i_t(M)$ is canonically placed
in $(\tm,\tg)$ as a totally geodesic submanifold for all values of $t$, 
and $i_t(M)$ is  disjoint from $i_{t'}(M)$ for all $t\neq t'$. 
Circles can be canonically placed into space forms with the placing 
corresponding to different critical values, so generally speaking, the
same manifold can admit different canonical placings into the ambient 
background.

\subsection{The mean curvature}
We let $\Psi(M)$ be the integral of the pointwise squared norm of
$H$:
\begin{equation}
\Psi(M)= \int_{M} \| H \|^2 \, d\mu \, . \label{mc}
\end{equation}
This functional is generally named after Willmore, 
who used it to study surfaces in $\mb{R}^3$ \cite{wi}. 
Such a functional was studied earlier by Blaschke \cite{bla}, and also by
Thomsen \cite{tho}.

If $f: (-a,a)\times M \rightarrow \tm$ is a family of 
isometric deformations of the immersion as above, we now compute the 
$t$-derivative of 
$\Psi(M_t)= {\displaystyle \int_{M_t} \| H_t\|^2\, d\mu_t}$ at $t=0$.
Since the ideas are similar to those used above in order to differentiate
(\ref{sf}), we shall be very brief this time around.

The key to the calculation is given by differentiating along the integral
curves of $T$ to obtain
$$
\begin{array}{rcl}
{\displaystyle \frac{d}{dt}\< H,H\>\mid_{t=0}} & = & 
2\< \n_T^{\tg} H,H\>+2\dot{g}^{ij}\< \a (e_i, e_j), H)\> \\ 
& = & 
2\< \n_T^{\tg} ( \n^{\tg}_{e_i}{e_i} -\n^{g}_{e_i}{e_i}) ,H\>
+2\dot{g}^{ij}\< \a (e_i, e_j),H\> \, .
\end{array} 
$$
Then we have the following:

\begin{theorem} \label{ele2}
Let $f: (-a,a)\times M \rightarrow \tm$ be a deformation of 
an isometrically immersed submanifold $f: M\rightarrow f(M)\hookrightarrow 
\tm$ into the Riemannian manifold $(\tm, \tg)$. We set $M_t=f(t,M)$, have
$M_0=M$, and let $\{ e_1, \ldots , e_n\}$ be an orthonormal frame of $M_t$ 
for all $t\in (-a,a)$. Then the infinitesimal variation of {\rm (\ref{mc})} is 
given by
$$
\begin{array}{rcl}
{\displaystyle \frac{d \Psi(M_t)}{dt}\mid_{t=0}} & = & 
{\displaystyle \int_{M} \< 2\n_{e_i}^{\tg}\n_{e_i}^{\tg}H+
2R^{\tg}(H,e_i)e_i
 -\| H \|^2 H,T\> d\mu } \vspace{1mm} \\
& & + {\displaystyle \int_{M} \< 4\< \a(e_i,e_j),H\>
\a(e_i,e_j)-\n \| H \|^2 , T\>d\mu 
} \vspace{1mm} \\
& & -2{\displaystyle \int_{M} (e_i\< T,e_j\> +e_j\<T,e_i\>)
\< \a (e_i, e_j),H\>d\mu }
\, .
\end{array}
$$
The isometrically immersed submanifold 
$f: M\rightarrow f(M)\hookrightarrow \tm$ satisfies the Euler-Lagrange 
equation of {\rm (\ref{mc})} if, and only if, the normal
component of the vector 
\begin{equation}
2\n_{e_i}^{\tg}\n_{e_i}^{\tg}H +2R^{\tg}(H,e_i)e_i
 - \| H \|^2 H + 4\< \a(e_i,e_j),H\>
\a(e_i,e_j)
\label{el2}
\end{equation}
is identically zero.
\qed
\end{theorem}

\subsection{Hypersurfaces}
For the case of a hypersurface, we proceed to re-express the Euler-Lagrange
equation of Theorems \ref{ele} \& \ref{ele2} as a scalar partial 
differential equation on the mean curvature function. When $(\tm,\tg)$ is 
assumed to be Einstein, this scalar equation can be expressed fully in terms
of the principal curvatures and the sectional curvatures of all normal
sections in the frame.

We write the mean curvature vector in the form $H=h\nu $, for some scalar
function $h$ and normal vector $\nu$. If $A$ is the shape operator, 
Wiengarten's identity (\ref{we}) implies that 
$\< \a(e_i,e_j),N\>=\< A_N(e_i),e_j\>$. By Codazzi's equation (\ref{co}),
we obtain that  
\begin{equation} \label{pr}
\begin{array}{rcl}
\< \n^{\tg}_{e_i}\n^{\tg}_{e_j}\a(e_i,e_j),\nu\>  & = & 
\< \n_{e_i}^{\tg}( (R^{\tg}(e_j,e_i)e_j)^{\nu}+\n^{\tg}_{e_i}\a
(e_j,e_j)+A_{\a(e_j,e_j)}e_i),\nu\> \\ & & \mbox{}\hfill  
-\< \n_{e_i}^{\tg}A_{\a(e_i,e_j)}e_j,\nu\>  \\
& = & 
 \< \n_{e_i}^{\tg}(R^{\tg}(e_j,e_i)e_j)^{\nu}),\nu\> +\< \n^{\tg}_{e_i}
\n^{\tg}_{e_i}
h \nu , \nu \> + \< \n_{e_i}^{\tg}A_{H}e_i,\nu\> \\ & & 
\mbox{}\hfill -\< \n_{e_i}^{\tg} A_{\a(e_i,e_j)}e_j ,\nu\> \, .
\end{array} 
\end{equation}

The various terms in this expression can be further developed. We have that 
$$
\begin{array}{rcl}
\< \n^{\tg}_{e_i} \n^{\tg}_{e_i} h \nu , \nu \> & = & 
\< \n_{e_i}^{\tg}( e_i(h)\nu+ h\n_{e_i}^{\tg}\nu), \nu \> \\
& = & \< ({e_i}(e_i (h)))\nu +2e_i(h) \n^{\tg}_{e_i} \nu + h \n_{e_i}^{\tg}
\n_{e_i}^{\tg}\nu, \nu \> \\ & = & -\Delta h - h \| \n_{e_i}^{\tg}\nu \|^2 \, ,
\end{array}
$$
and the Weingarten identity (\ref{we}) yields that
\begin{equation}
\begin{array}{rcl}
\< \n^{\tg}_{e_i} \n^{\tg}_{e_i} h \nu , \nu \> & = & 
-\Delta h - h \| A_{\nu}({e_i})\|^2 -h\| \n_{e_i}^{\nu}\nu \|^2 \\
& = & -\Delta h - h\, {\rm trace} \, A_{\nu}^2 \, ,
\end{array} \label{key}
\end{equation}
because $2\< \n_{e_i}^\nu \nu , \nu\>=0$, and so 
$\n_{e_i}^\nu \nu $ is the zero vector.

We also have that
$$
\< \n_{e_i}^{\tg }A_{H}(e_i), \nu \> = \< \n_{e_i}^{\tg}(h A_{\nu}(e_i)),
\nu\>= - h\< A_{\nu}(e_i), \n_{e_i}^{\tg}\nu \> = h\, 
{\rm trace}\, A_{\nu}^2 \, .
$$
Finally, 
$$
\< \n^{\tg}_{e_i}A_{\a(e_i,e_j)}e_j,\nu\> =-
 \<A_{\a(e_i,e_j)}e_j,\n_{e_i}^{\tg}
\nu\> = \< A_{\a(e_i,e_j)}e_j,A_{\nu}(e_i)\> \, .
$$
We now set $\a(e_i,e_j)=h_{ij}\nu$, and choose the frame so that $A_\nu$ is
diagonal to see that
$$\< \n^{\tg}_{e_i}A_{\a(e_i,e_j)}e_j,\nu\> =h_{ij}\<A_{\nu}e_j,A_{\nu}e_i\> =
 {\rm trace}\, A^3_{\nu } \, .
$$
Using these results in (\ref{pr}), we obtain that
\begin{equation}
\< \n^{\tg}_{e_i}\n^{\tg}_{e_j}\a(e_i,e_j),\nu\>  =  
 \< \n_{e_i}^{\tg}(R^{\tg}(e_j,e_i)e_j)^{\nu}),\nu\> -\Delta h 
-{\rm trace}\, A_{\nu}^3  \, . \label{lu1}
\end{equation}

Once again, by choosing a frame that diagonalizes 
$A_{\nu}$, we can easily see that
\begin{equation}
\< \a(e_i,e_j),\a(e_l,e_j)\>\<\a (e_l,e_i),\nu\>= {\rm trace}\, A_{\nu}^3\, . 
\label{ul1}
\end{equation}

We now have the following results. The first states the Euler-Lagrange 
equation of the functionals (\ref{sf}) and (\ref{mc}) for hypersurfaces in 
the case when $(\tilde{M},\tg)$ 
is Einstein. We then draw some consequences. 

\begin{theorem} \label{th1}
Let $M$ be a hypersurface of an Einstein manifold $(\tm,\tg)$. Assume that
$k_1, \ldots , k_n$ are the principal curvatures, with associated orthonormal 
frame of principal directions $e_1, \ldots, e_n $. Let $\nu$ be a normal
field along $M$. Then $M$ is a critical point of the 
functional {\rm (\ref{sf})} if, and only if, 
$$
2\Delta h = 2 (k_1 K_{\tg}(e_1,\nu)+\cdots + k_n K_{\tg}(e_n,\nu))
 -h\| \a\|^2 +
2(k_1^3 + \cdots + k_n^3) \, ,
$$
and a critical point of the functional {\rm (\ref{mc})} if, and only if, 
$$
2\Delta h = 2h(K_{\tg}(e_1,\nu)+\cdots + K_{\tg}(e_n,\nu))
 +2h\| \a\|^2 -h^3 \, .
$$
\end{theorem}

Notice that we have $h=k_1 + \cdots + k_n$ and $\| \a\|^2 = k_1^2 + 
\cdots + k_n^2$, respectively.

{\it Proof}. We show that the vanishing of the normal component of (\ref{el}) 
yields precisely the first stated critical point equation.  

The term $\<R^{\tg}(\a (e_i,e_j),e_j)e_i$ gives rise to sectional curvatures.
Indeed, since $\a (e_i,e_j)=k_i \delta_{ij}\nu$, we have that 
\begin{equation}
\<R^{\tg}(\a (e_i,e_j),e_j)e_i ,\nu \>=k_i \<R(\nu, e_i)e_i, \nu\> =
k_i K_{\tg}(e_i,\nu)\, . \label{ul2}
\end{equation}

The first term in the right side of (\ref{lu1}) is zero.
Indeed,
$\<\n_{e_i}^{\tg}(R^{\tg}(e_j,e_i),e_j)^{\nu} ,\nu \>= 
e_i\<(R^{\tg}(e_j,e_i),e_j)^{\nu} ,\nu \>-\<(R^{\tg}(e_j,e_i),e_j)^
{\nu} ,\n_{e_i}^{\tg}\nu \>=e_i\<(R^{\tg}(e_j,e_i),e_j)^{\nu} ,\nu \>$,
because $\n_{e_i}^{\tg}\nu$ is tangent to $M$, and since the 
metric $\tg$ is Einstein, if $s_{\tg}$ is the scalar curvature, we have 
that 
$$
\<(R^{\tg}(e_j,e_i),e_j)^{\nu} ,\nu \>=
-\<R^{\tg}(e_j,e_i)\nu, e_j\>=-\frac{s_{\tg}}{\tn} \tg( e_i,\nu)=0\, .
$$

The desired equation results by using (\ref{lu1}), (\ref{ul1}) and
(\ref{ul2}) in the computation of the normal
component of (\ref{el}).

The vanishing of the normal component of (\ref{el2}) yields the stated critical
point equation for (\ref{mc}). In order to derive this, we begin by 
observing that the normal component of $R^{\tg}(H,e_i)e_i$ leads to
$h\sum_i K_{\tg}(e_i, \nu)$. The highest order derivative
term arises directly by (\ref{key}). 
The term $\< \a(e_i,e_j),H\>\<\a(e_i,e_j),\nu\>$ yields $h\| \alpha \|^2$, 
which simplifies with the lower order cubic term in the right side of 
(\ref{key}). 
\qed

We denote by $S_c^n$ the $n$th dimensional space form of curvature 
$c$. 

\begin{corollary} \label{c4}
Let $k_1, \ldots , k_n$ be the principal curvatures of a hypersurface $M$ in
$S_c^{n+1}$. Then $M$ is a critical point for the functional {\rm (\ref{sf})} 
if, and only if, its mean curvature function $h$ satisfies the equation
$$
2\Delta h = 2 ch -h\| \a\|^2 +
2(k_1^3 + \cdots + k_n^3) \, ,
$$
while $M$ is a critical point of the functional {\rm (\ref{mc})} if, and
only if, its mean curvature function $h$ satisfies the equation
$$
2\Delta h = 2cn h +2h\| \a\|^2 -h^3 
$$
instead.
\qed
\end{corollary}

\begin{corollary}
Austere hypersurfaces of a space form $S_c^n$ are critical
points of the functional {\rm (\ref{sf})}, and so they are canonically placed
into $S_c^n$. They are absolute minima for
the functional {\rm (\ref{mc})}, functional that is also minimized by any
other minimal hypersurface.
\qed
\end{corollary}

Austere manifolds were introduced by Harvey and Lawson \cite{hala}
in their study of Calibrated Geometries. They are submanifolds of
a Riemannian manifold for which the eigenvalues of its fundamental
form in any normal direction occur in opposite signed pairs.


\begin{remark}
Let $S_r$ be the $r$th symmetric function of the 
principal curvatures of the hypersurface $M$. For a given real valued function 
$f$ in $n$ variables, Reilly studied the variation of 
$$
\int_M f(S_1,\ldots, S_n)d\mu 
$$ 
under deformations \cite{rcr}. He deduced the variational formula for the 
$L^2$-norm of the second fundamental form for hypersurfaces in space forms,
obtaining the equation given in Corollary \ref{c4}.
Reilly applied this result in the case of ${\mathbb S}^3$, where via 
the identity $a^3+b^3=(a+b)(a^2-ab+b^2)$, he concluded that the only
critical points on the $3$-sphere whose principal curvatures never change
sign are minimal surfaces. This is of course included in 
our conclusion above, as in dimensions one and two, the notions of austere
and minimal submanifolds coincide. We will see next 
that even if we require that all principal curvatures be constant functions,
there are critical hypersurfaces of (\ref{sf}) on the sphere that are not
minimal.
\qed
\end{remark}

\begin{example}\label{ex8}
Even in codimension one, the class of critical points of (\ref{sf}) is 
larger than the class of austere hypersurfaces. In order to see that, we 
use the family of 
isoparametric hypersurfaces of Nomizu \cite{no}, given by
$M_t^{2n}=\{z \in {\mathbb S}^{2n+1}: F(z)={\rm cos}^2 2t\}$
for $0<t<\pi/4$, where
$F(z)=(\| x\|^2-\| y\|^2)^2 +4\< x, y\>^2$, $z=x+iy$, $x,y 
\in {\mathbb R}^{n+1}$. For each $t$, this hypersurface has principal 
curvatures $k_1=(1+\sin{2t})/\cos{2t}$ and $k_2=(-1+\sin{2t})/\cos{2t}$, with
multiplicity $1$ each, and $k_3=\tan{t}$ and $k_4=-\cot{t}$ with multiplicity 
$n-1$ each, respectively. We thus have that $h=h(t)= 
k_1(t)+k_2(t)+(n-1)(k_3(t)+k_4(t))$.

When $n=2$, $M_t^4$ is an austere hypersurface
for the value of $t$ given by the root of 
$h(t)=k_1(t)+k_2(t)+k_3(t)+k_4(t)$ in the interval $(0,\pi/4)$, which
coincides with the root of 
$2h(t)-\| \a\|^2 h(t) +2(k_1^3(t)+k_2^3(t)+k_3^3(t)+k_4^3(t))$, and which
is close to $t=0.3926990817$. The corresponding critical value of
the functional (\ref{sf}) is $12$ times the volume of the manifold. 

However, for $n\geq 3$, the functions 
$2h(t)-\| \a\|^2 h(t) + 2(k_1^3(t)+k_2^3(t)+(n-1)(k_3^3(t)+k_4^3(t)))$ and
$h(t)$ do have unique distinct roots on the interval $(0,\pi/4)$.
Thus, we have a critical point of $(\ref{sf})$ that is not a minimal
hypersurface. For instance, when $n=3$, the root of the first function is
$t=0.3775786497$, and the corresponding critical value of (\ref{sf}) is
$18.57333958$ times the volume. In this case, the function $h(t)$ has a root at
$t=0.4776583091$.

When looking at the functional (\ref{mc}) instead, we can use these 
isoparametric hypersurfaces also to show that the class of critical points 
is larger than the class of
minimal hypersurfaces. Indeed, for $n=4$ the function 
$2n+ 2\| \a\|^2 -h^2(t)$ has a root at $t=0.2153460562$, where
$h(t)$ is not zero. Thus, this $M_{t}^{8}$ is a critical hypersurface of
(\ref{mc}) in ${\mathbb S}^9$. Its critical value is $147.3776409$ times
the volume. Let us notice that this function is strictly positive for $n$ equal
$2$ and $3$, respectively.
\qed
\end{example}

\begin{remark}
In \cite{si}, Simons computed the Laplacian of $\a$ for a minimal 
variety of dimension $p$ immersed into an Euclidean $n$-sphere. Using 
the nonnegativity of the Laplacian, he concluded
that for a minimal $p$-variety immersed into the $n$-sphere such that
$0\leq \| \a\|^2 \leq p/q$, $q=2-1/(n-p)$, then either $\| \a\|^2=0$ or 
$\| \a \|^2=p/q$, respectively. 
The Nomizu austere hypersurface $M^4_t\subset \mb{S}^5$ above
has $\| \a\|=12$, a constant that falls 
outside the range covered by Simons' result. 
This one of the isoparametric families 
of Nomizu appears as one of the examples studied by Cartan \cite{cart}. 
In \cite{pete}, Peng and Terng showed that among all of the 
Cartan's isoparametric families, the only possible values of $\| \a \|^2$
are $0$, $n-1$, $2(n-1)$, $3(n-1)$ or $5(n-1)$, respectively.  Their 
result produces examples of some rather large values of constant 
$\|\a\|^2$ for 
minimal hypersurfaces of the $n$-sphere. In Example \ref{ex8} we had seen
already a rather large value of constant $\| \a \|^2$ for a critical 
hypersurface of (\ref{sf}) in $\mb{S}^{2n+1}$, $n\geq 3$.
\qed
\end{remark}

\begin{example}\label{ex7}
Since the orbit of a compact connected group of isometries whose volume is
extremal among nearby orbits of the same type is a minimal submanifold
\cite{hs}, studying the diagonal action of ${\mathbb S}{\mathbb O}(3,
{\mathbb R})$ in ${\mathbb R}^3 \times  {\mathbb R}^3$, we may conclude 
\cite{bl} that
${\mathbb S}^2(\sqrt{1/2})\times {\mathbb S}^2(\sqrt{1/2})\subset
{\mathbb S}^5$ is a minimal hypersurface. Here, ${\mathbb S}^n(r)$ denotes
the $n$-sphere of radius $r$, with ${\mathbb S}^n(1)={\mathbb S}^n$. 
The principal curvatures of the said minimal
submanifold are
$1$, $1$, $-1$, and $-1$, respectively. Thus, this manifold is austere, and
consequently, a critical point of (\ref{sf}).

In fact, more is true: of all the products 
${\mathbb S}^m(\sqrt{n/(m+n)})\times {\mathbb S}^{n}(\sqrt{m/(m+n)})\subset
{\mathbb S}^{m+n+1}$, the only ones that are critical points of
(\ref{sf}) must have $m=n$, in which case, they are austere. This follows
easily from the explicit calculation of the principal curvatures, which
turn out to be $k_1= \ldots =k_m=\sqrt{(m/n)}$ and
$k_{m+1}= \ldots =k_{m+n}=-\sqrt{(n/m)}$, respectively. 
Cf. \cite{glw}, \S 3, and \cite{hl}, p. 366. All of the alluded to product of 
spheres are critical points of the functional 
${\displaystyle \int_M \left( \| \a\|^2- \frac{1}{n}\| H\|^2
\right)^{\frac{m+n}{2}} d\mu}$,
termed Willmore in these references. The functional (\ref{sf}) distinguishes 
the symmetric Clifford minimal tori out of these products. These special tori
are Einstein manifolds \cite{rss}.  
\qed
\end{example}

\subsection{Submanifolds of arbitrary codimension}
We now consider submanifolds of arbitrary codimension $q$, and re-express the 
Euler-Lagrange equation of Theorems \ref{ele} \& \ref{ele2} 
as $q$-scalar equations. 

For that, we choose an orthonormal frame $\{ \nu_1,\nu_2, \ldots, \nu_q\}$ 
for the normal bundle $\nu(M)$ so that the mean curvature vector is 
parallel to $\nu_1$. 

The higher codimension analogue of (\ref{lu1}) along $\nu_1$ is obtained 
once again by Codazzi's 
equation (\ref{co}), and is given by 
\begin{equation}
\< \n^{\tg}_{e_i}\n^{\tg}_{e_j}\a(e_i,e_j),\nu_1\>  =  
 \< \n_{e_i}^{\tg}(R^{\tg}(e_j,e_i)e_j)^{\nu}),\nu_1\> -\Delta h - h
\| \n^{\nu}_{e_i}\nu_1 \|^2 -{\rm trace}\, A_{\nu_1}A_{\nu_k}^2 
\, \label{ul3}
\end{equation}
(with our summation convention on repeated indexes, the last trace in the
right side above is in effect the trace of $ A_{\nu_1}\sum_{k=1}^q 
A_{\nu_k}^2$). This time around, $\sum_i \n_{e_i}^{\nu}\nu_1$ is not
necessarily the zero vector. The other difference 
arises in manipulating the expression 
$\< A_{\a(e_i,e_j)}e_j,A_{\nu_1}(e_i)\>$. We now
decompose the normal vector $\a(e_i,e_j)$ as $\a(e_i,e_j)=h_{ij}^{k}\nu_k$
instead of the one dimensional expression along $\nu$ used earlier. By our
assumptions, $h_{ii}^1=h$ while $h_{ii}^k=0$ for $k\geq 2$, and 
the cubic trace term in (\ref{ul3}) results by choosing frames that
diagonalize $A_{\nu_k}$.

Similarly, for any $m$ in the range $1\leq m\leq q$, we have that
\begin{equation}
\< \a(e_i,e_j),\a(e_l,e_j)\>\<\a (e_l,e_i),\nu_m\>=
{\rm trace}\, A_{\nu_m}A_{\nu_k}^2\, , \label{bl1}
\end{equation}
and so the vanishing of the component of (\ref{el}) along $\nu_1$ yields 
the equation
\begin{equation}
\begin{array}{rcl}
0 & = & 2\<R^{\tg}(\a (e_i,e_j),e_j)e_i ,\nu_1 \>
+2\< \n_{e_i}^{\tg}(R^{\tg}(e_j,e_i)e_j)^{\nu}),\nu_1\>
 -2\Delta h 
-2h\| \n_{e_i}^{\nu}
\nu_1 \|^2 \\ & &  -
h\| \a\|^2 +2{\rm trace}\, A_{\nu_1}A_{\nu_k}^2 \, . 
\end{array}\label{eq1}
\end{equation}

The equation arising from the vanishing of the component of (\ref{el}) 
along $\nu_m$ for $m\geq 2$ follows by a similar argument but 
is slightly different, the reason being that the trace of $\a$ is 
orthogonal to the said direction. For $H=h\nu_1$, and
so 
$$
\< \n^{\tg}_{e_i} \n^{\tg}_{e_i} H , \nu_m \> =2e_i(h)
\< \n_{e_i}^{\nu}\nu_1,\nu_m\> + he_i 
\< \n_{e_i}^{\nu}\nu_1,\nu_m\>
-h\, {\rm trace}A_{\nu_1}A_{\nu_m}-h\< \n_{e_i}^{\nu}
\nu_1, \n_{e_i}^{\nu} \nu_m \> \, .
$$

We also have that
$$
\< \n_{e_i}^{\tg}A_H (e_i),\nu_m\> =h\< A_{\nu_1}(e_i),A_{\nu_m}(e_i)\>
=h \, {\rm trace} A_{\nu_1}A_{\nu_m} \, .
$$
And since 
$$
\< \n^{\tg}_{e_i}A_{\a(e_i,e_j)}e_j,\nu_m\> =-
 \<A_{\a(e_i,e_j)}e_j,\n_{e_i}^{\tg}
\nu_m\>  =  \< A_{\a(e_i,e_j)}e_j,A_{\nu_m}(e_i)\>\, ,
$$
choosing frames that diagonalize the $A_{\nu_k}$s, we obtain
that
$$
\< \n^{\tg}_{e_i}A_{\a(e_i,e_j)}e_j,\nu_m\>={\rm trace} A_{\nu_m}A_{\nu_k}^2
\, .
$$
Thus, the vanishing of the component of (\ref{el}) along $\nu_m$ yields 
the equation
\begin{equation}
\begin{array}{rcl}
0 & = & 2\< R^{\tg}(\a (e_i,e_j),e_j)e_i , \nu_m\>+2
\<\n_{e_i}^{\tg}(R^{\tg}(e_j,e_i)e_j)^{\nu}), \nu_m \> 
+ 4e_i(h)
\< \n_{e_i}^{\nu}\nu_1,\nu_m\> \\ & & + 2he_i 
\< \n_{e_i}^{\nu}\nu_1,\nu_m\>-2h\< \n_{e_i}^{\nu}
\nu_1, \n_{e_i}^{\nu} \nu_m \> 
+2{\rm trace}\, A_{\nu_m}A_{\nu_k}^2 \, .
\end{array}
\label{lu3}
\end{equation}

The vanishing of the various normal components of (\ref{el2}) 
follow by similar considerations. We simply state the final result.
When $q=1$, the statement that follows displays also the Euler-Lagrange 
equation of (\ref{sf}) and (\ref{mc}), respectively, in full generality,
without the Einstein assumption on $(\tm,\tg)$ used in Theorem \ref{th1}.

\begin{theorem} \label{th2}
Let $M$ be an $n$-manifold of codimension $q$ isometrically immersed 
in $(\tm,\tg)$. Let $\{ e_1, \ldots, e_n\}$ be an orthonormal frame of
$M$. Suppose that $\{ \nu_1,\ldots, \nu_q\}$
is an orthonormal frame of the normal bundle of the immersion 
such that $H=h\nu_1$. Then, $M$ is a critical point of the 
functional {\rm (\ref{sf})} if, and only if, 
$$
2\Delta h  =  2\<R^{\tg}(\a (e_i,e_j),e_j)e_i +
\n_{e_i}^{\tg}(R^{\tg}(e_j,e_i)e_j)^{\nu},\nu_1 \> -2h\| \n_{e_i}^{\nu}
\nu_1 \|^2 -
h\| \a\|^2 +2{\rm trace}\, A_{\nu_1}A_{\nu_k}^2 \, ,
$$
and for all $m$ in the range $2\leq m \leq q$, we have that
$$
\begin{array}{rcl}
0 & = & 2\<  R^{\tg}(\a (e_i,e_j),e_j)e_i +
\n_{e_i}^{\tg}(R^{\tg}(e_j,e_i)e_j)^{\nu}+2e_i(h)\n_{e_i}^{\nu}\nu_1,\nu_m\>
\\ & &  +2he_i \< \n_{e_i}^{\nu}\nu_1,\nu_m\> 
- 2h\< \n_{e_i}^{\nu}
\nu_1, \n_{e_i}^{\nu} \nu_m \>  
+2{\rm trace}\, A_{\nu_m}A_{\nu_k}^2\, .
\end{array}
$$ 
$M$ is a a critical point of the 
functional {\rm (\ref{mc})} if, and only if, 
$$
2\Delta h  =  2h(K_{\tg}(e_1,\nu_1) + \cdots +K_{\tg}(e_n,\nu_1) )
-2h\| \n_{e_i}^{\nu}\nu_1 \|^2 -
h^3 +2h\, {\rm trace}\, A_{\nu_1}^2 \, ,
$$
and for all $m$ in the range $2\leq m \leq q$, we have that
$$
\begin{array}{rcl}
0 & = & 2h\< R^{\tg}(\nu_1,e_i)e_i,\nu_m\>+4e_i(h)\< \n_{e_i}^{\nu}\nu_1,
\nu_m\>+2he_i \< \n_{e_i}^{\nu}\nu_1,\nu_m\>
 \\ & & - 2h\< \n_{e_i}^{\nu}
\nu_1, \n_{e_i}^{\nu} \nu_m \> 
+2h\, {\rm trace}\, A_{\nu_1}A_{\nu_m}\, .
\end{array}  
$$
\qed
\end{theorem}

\begin{remark}
For $q=1$ above, and setting $\nu=\nu_1$, we have that 
$$
\< \n_{e_i}^{\tg}(R^{\tg}(e_j,e_i)e_j)^{\nu},
\nu\>=-e_i (r_{\tg}(e_i,\nu))\, .
$$
If $(\tm,\tg)$ were to be Einstein, cf. with Theorem \ref{th1}, 
this expression would be identically zero. 
\qed
\end{remark}
\medskip

\subsection{The extrinsic scalar curvature}
We let $\Theta(M)$ be given by the total extrinsic scalar curvature of the 
immersion. This is the remaining piece in our analysis of the integral of the various
summands in the right side of (\ref{sc}):
\begin{equation}
\Theta(M)= \int_{M} \sum_{1\leq i,j \leq n} K_{\tg} 
(e_i,e_j) \, d\mu \, . \label{esc}
\end{equation}

\begin{theorem}
Let $M$ be an $n$-manifold of codimension $q$ isometrically immersed 
in $(\tm,\tg)$. Let $\{ e_1, \ldots, e_n\}$ be an orthonormal frame of
$M$. Suppose that $\{ \nu_1,\ldots, \nu_q\}$
is an orthonormal frame of the normal bundle of the immersion 
such that $H=h\nu_1$. Then, $M$ is a critical point of the 
functional {\rm (\ref{esc})} if, and only if, 
$$
h\sum_{1\leq i,j\leq n}K_{\tg}(e_i,e_j)=0\, .
$$
\end{theorem}

{\it Proof}. Let $T$ be the variational vector field of the deformation, and
$\{ e_1, \ldots, e_n\}$ a normal frame of $M_t$ at $(p,t)$  
for sufficiently small $t$, as used above. Then the value of  
$K_{\tg}(\dot{e}_i,e_j)+K_{\tg}(e_i,\dot{e}_j)$ at $t=0$ is zero because
the $e_j$ component of $\dot{e}_i \mid_{t=0}$ is minus the $e_i$ component
of $\dot{e}_j\mid_{t=0}$. The result follows using this 
fact and (\ref{vf}) in the differentiation of (\ref{esc}).
\qed


\begin{corollary} \label{c313}
Let $(M,g)$ be a manifold isometrically embedded into $(\tm, \tg)$. 
If the extrinsic scalar curvature of $M$ is nowhere zero, $M$ is a 
critical point of {\rm (\ref{esc})} if, and only if, it is minimal. 
In particular, if $\tg$ has positive or negative sectional curvature, 
$M$ is a critical point of ${\rm (\ref{esc})}$ if, and only if, $M$ is minimal.
\end{corollary}

The condition on the sectional curvature of $\tg$ above can be relaxed somewhat
while maintaining the same conclusion. We need to ask that the set of points
carrying planar sections be ``small'' in a sense we leave to the reader to formulate.

For isometric immersions of $M$ into $(\tm,\tg)$, we may 
consider the functional
\begin{equation}
M \mapsto {\mathcal S}(M)=\int_{M}(\| \a\|^2 -\| H\|^2) 
d\mu \, , 
\label{ssf}
\end{equation}
Its critical point equation may be obtained by taking the difference of the
critical point equations for the functionals $\Pi(M)$ and $\Psi(M)$ given by
Theorem \ref{th2} above. We shall use this functional later on to
show the existence of bulk backgrounds with 
isotopic critical submanifolds of different critical values. These examples
exhibit the same flavour of some alternatives theories to 
compactification currently of interest to physicists \cite{rasu}.

\subsection{Canonical homology representatives}
Assuming that the topology of $\tm$ is nontrivial, 
we revisit an old question of Steenrod inquiring about  
the shape that a homology class has. This paraphrases
the formulation of Steenrod's question in the language of 
D. Sullivan \cite{dennis}.

We assume that $\tm$ is a smooth connected oriented manifold, and let
$D$ be a class in $H_n(\tm;{\mathbb Z})$. A closed oriented manifold 
$M=M^n$ is said to
represent
$D$ \cite{th} if there exists an embedding $i: M \hookrightarrow \tm$ such that
$i_{*}[M]=D$.

We choose and fix a homology class
$D\in H_n(\tm ;{\mathbb Z})$ in $\tm$,
and define the space
\begin{equation}
{\mathcal M}_D(\tm)=
\{M:  M \hookrightarrow \tm \;
\text{is an isometric embedding,}\; [M]=D\} \, . \label{hcs}
\end{equation}
Since the metric $\tg$ on $\tm$ is fixed, we drop it altogether
from the notation here.

Points in $\mc{M}_D(\tm)$ are manifolds that can have different
topology. The stability of the class of smooth embeddings of 
a compact manifold into a fixed background allows us to introduce a topology 
on $\mc{M}_D(\tm)$ by topologizing the small deformations of a point
in $\mc{M}_D(\tm)$ to define neighborhoods of the point in question.

\begin{remark}\label{thre}
A celebrated result of Thom \cite{th} shows that for
any integral homology class
$D$ there exists a nonzero integer $N$ such that the class
$N\cdot D$ is represented by a submanifold $M$. We assume here 
that the class $D$ is such that $N$ can be chosen to be $1$,
and so the space ${\mathcal M}_D(\tm)$ is nonempty.
This nontrivial issue for integer homology classes
``evaporates'' when using field coefficients; for
Thom derives from the result above \cite[Corollary II.30]{th}  
that the homology groups with rational or real coefficients admit 
bases of systems of elements represented by submanifolds. 
He proved that any class in $H_{i}(\tm; {\mathbb Z})$, $i\leq 6$, is 
realizable by a submanifold \cite[Theorem II.27, footnote 9
p. 56]{th}, and then deduced that for oriented 
manifolds of dimension at most $9$, any integral class can be so realized
\cite[Corollary II.28, footnote 9]{th}. Thom exhibited an example of a 7 
dimensional integer class in a manifold of dimension $14$ 
that is not realizable by a submanifold. Optimal examples of nonrealizable
$7$-classes on manifolds of dimension $10$ had been given in \cite{bhk}.
See Remark \ref{bhk} below. 

It is simple to see that in codimension $q=\tn -n=1$ or $2$, 
there is no obstruction to 
representing $D\in H_{n}(\tm; {\mathbb Z})$ by an embedded submanifold, and
${\mathcal M}_D(\tm)$ will therefore be always nonempty. For instance, if 
$q=2$, the Poincar\'e dual of $D$ will be an element of $H^2(\tm;
{\mathbb Z})$, group that is in 1-to-1 correspondence with the group of
isomorphism classes of complex line bundles over $\tm$ via the first 
Chern class
mapping. The zero set of a generic section of the line bundle
that corresponds to the said Poincar\'e dual is a
smooth submanifold that represents $D$. The structure group of its normal
bundle is $\mb{S}\mb{O}(2)\cong \mb{U}(1)$.
\qed
\end{remark}
\medskip

The {\it nicest} metric element in the space ${\mathcal M}_D(\tm)$ in 
(\ref{hcs}) should describe the way the {\it homology class $D$ looks 
like} \cite{dennis}. This problem, thoroughly 
discussed by Thom via topological techniques, is now  
revisited by using a variational principle and provided with 
an alternative geometric answer. We have in mind the
definition of a functional
on ${\mathcal M}_D(\tm)$ whose minimum ---when it can be achieved---
defines the nicest element in question. Of course, such an optimization problem
will depend upon the functional we choose to use for the purpose.

If $D$ is a one dimensional homology class, its natural optimal representative
is the shortest ambient space geodesic loop in the class, a submanifold that when
provided with the induced metric has second fundamental form that 
is identically zero, and so it yields an absolute minimizer of the functional 
$\Pi$ in (\ref{sf}) defined over the free homotopy class of $D$.  
(Notice that this is also the natural answer to the
question of finding optimal representatives of elements of the fundamental 
group instead.) However, even in this simple case, the functional $\Pi$ may have other 
critical points of higher critical value, and these nonminimizer curves would encode 
additional information on the curvature of the metric in the background manifold.  
In higher dimension, on the other hand, the curvature itself of the
isometrically immersed $(M,g)$ 
should play more of a role in determining the best representative of $D$, even if the
latter were to be defined as a minimizer of $\Pi$ within the class. 

The curvature functionals to optimize in order to find canonical 
representatives of $D$ are plentiful, but the identity (\ref{sc}) suggests a natural choice
to make. We start by restricting the domain of (\ref{sf}), and seek critical points of
\begin{equation}
\begin{array}{ccl}
{\mathcal M}_{D}(\tm ) & \rightarrow & {\mathbb R} \\
M & \mapsto & {\displaystyle \mc{T}(M)=\int_{M} \tg(\a,\a) d\mu_g }
\end{array}\, . \label{cpp}
\end{equation}
Since the integrand is bounded below and often $D$ has volume minimizer representatives,
it is natural to seek critical points of (\ref{cpp}) that 
are minima and have the smallest possible volume among these. Critical points of (\ref{cpp}) could 
have some connected components that bound (which would 
produce homologically trivial critical points of (\ref{sf}) in their own 
right), or connected components that cancel each other out partially in 
homology to produce a representative of $D$. It is then natural to impose the mentioned 
additional constraint, and seek minimizers of $\Pi$ in $\mc{M}_D(\tm)$ that have 
the smallest volume. A minimum $M$ of smallest volume,
should it exist, will be said to be a {\it canonical representative} of $D$. 

Notice that if a canonical representative $M$ of the class $D$ 
were to be a volume minimizer also, (\ref{sc}) would imply that the metric on $M$   
is as close as it can be to an Einstein metric among all the metrics on $M$ induced by
the ambient space metric $\tg$ on $\tm$.  

\begin{remark}
It might be natural to restrict the domain of the functional (\ref{cpp}) to an open 
subset of $\mc{M}_D(\tm)$ according to the topology of its elements, open subset that
encodes additional information about differential invariants of the 
background manifold $\tm$. In that case, we might want to redefine a canonical 
representative of $D$ accordingly. 
\qed
\end{remark}

\subsection{An example: the fibration ${\mathbb S}^3 
\hookrightarrow Sp(2)\stackrel
{\pi_{\circ}}{\rightarrow} {\mathbb S}^7$} \label{s37}
If ${\mathbb H}=\{q=q_0 +q_1i + q_2 j + q_3k: \; (q_0,q_1,q_2,q_3)\in 
{\mathbb R}^4\}$ is the (skew) field of quaternions, the
symplectic group $Sp(n)$ is the group of $n\times n$ matrices $M$ with 
entries in ${\mathbb H}$ such that
$\< Mx,My\>_{\mathbb H}=\< x,y\>_{\mathbb H}$. Here, 
$\< x,y\>_{\mathbb H}=\sum_i \bar{x}_i y_i $
is the standard quaternionic bilinear form on ${\mathbb H}^n$. 
$Sp(n)$ is a simple compact Lie group of dimension $2n^2+n$.

By identifying a quaternion with its corresponding vector of coordinates in
${\mathbb R}^4$, we see that $Sp(1)={\mathbb S}^3$.
On the other hand, 
\begin{equation}
Sp(2)=\{ \left( q_{ij}\right) \in {\mathfrak M}_2({\mathbb H}): \;
\sum_{l=1}^2 \bar{q}_{l i}q_{l j}=\delta_{ij }\, , 
\; 1\leq  i, j \leq 2 \} \, ,\label{eq3}
\end{equation}
and its Lie algebra is given by
\begin{equation}
{\mathfrak s}{\mathfrak p}(2)=\left\{ \left( q_{ij}\right) \in 
{\mathfrak M}_2(
{\mathbb H}): \; \bar{q}_{21}=-q_{12}\, , \; \bar{q}_{\a \a}=-q_{\a \a}\;
{\rm for} \; 1\leq  \a \leq 2 \right\} \, .\label{eq4}
\end{equation} 

We consider the free action of $Sp(1)$ on $Sp(2)$
defined by letting the quaternion $q\in Sp(1)$ act on $Sp(2)$ according 
to the rule
\begin{equation}
q \circ \left( q_{\a \b}\right) =
\left( \begin{array}{cc} q_{11} & q_{12} \\ q_{21} & q_{22} \end{array}\right)
\left( \begin{array}{cc} 1 & 0 \\ 0 & \bar{q} \end{array}\right) \, .
\label{ac0} 
\end{equation}
The principal fibration of $Sp(2)$ induced by the $\circ$-action 
has as base ${\mathbb S}^7$ with its standard differentiable
structure. The projection map $\pi_{\circ}$ of this fibration is 
given by
$$
\pi_{\circ}\left( \begin{array}{cc} q_{11} & q_{12} \\ q_{21} & q_{22} 
\end{array}\right) = \left( \begin{array}{c} q_{11} \\ q_{21}
                            \end{array}\right)
\, .
$$
The action of $Sp(2)$ on itself by left translations 
commutes with $\pi_{\circ}$, and so $Sp(2)$ acts in this manner by bundle 
morphisms of the fibration ${\mathbb S}^3 \hookrightarrow 
Sp(2)\stackrel
{\pi_{\circ}}{\rightarrow} {\mathbb S}^7$. We thus obtain an induced 
action of $Sp(2)$ on ${\mathbb S}^7$, which is defined by matrix 
multiplication when ${\mathbb S}^7$ is thought of as the subset of 
${\mathbb H}^2$ consisting of column vectors that have unit length.

Notice that if $\gamma_q :[0,1] \rightarrow \mb{S}^3$ is a smooth curve in 
$\mb{S}^3$ that begins at $q\in \mb{S}^3$ and ends at $1\in \mb{S}^3$, the fibration above can 
be made part of the 
$t$-parameter family of fibrations 
${\mathbb S}^3 \hookrightarrow 
Sp(2)\stackrel
{\pi_{{\circ}_t}}{\rightarrow} {\mathbb S}^7$ 
defined by the curve of $Sp(1)$-actions  
\begin{equation}
q \circ_t \left( q_{ij}\right) =
\left( \begin{array}{cc} q_{11} & q_{12} \\ q_{21} & q_{22} \end{array}\right)
\left( \begin{array}{cc} \overline{\gamma}_q(1-t) & 0 \\ 0 & 
\overline{\gamma}_q(t) 
\end{array}\right)\, , 
\label{ac1} 
\end{equation}
$t\in [0,1]$. The fibration corresponding to $t=1$ projects 
$(q_{ij})$ onto the point of $\mb{S}^7$ represented by its second column.
It a sort of rotated version of $\pi_{\circ}$.

The fibers and base of $\pi_{\circ}$ play special roles
in relation to the geometry defined by certain Riemannian metrics on $Sp(2)$.
We let $g_{\lambda}$ be the metric on the $7$-sphere obtained by scaling the 
vertical space of the Hopf fibration by the factor $\lambda$. Now, as 
a Lie group of dimension $10$, $Sp(2)$ has a $55$-dimensional family of
distinct left invariant metrics, and of these, there is a $2$-dimensional
family $g_{\lambda,\mu}$ that descends to $g_{\lambda}$ on ${\mathbb S}^7$
and that induces constant sectional curvature metrics on the 
$\mb{S}^3$-fibers of $\pi_{\circ}$. The metric $g_{\lambda, \mu}$ 
yields a Riemannian submersion
$(Sp(2), g_{\lambda,\mu}) \rightarrow ({\mathbb S}^7,g_{\lambda})$ whose
fibers are all isometric $3$-spheres.
For, since left invariant vector fields are uniquely defined by their 
values at the identity element of $Sp(2)$, it suffices to describe the
metric on ${\mathfrak s}{\mathfrak p}(2)$. Then $g_{\lambda,\mu}$ is defined by
\begin{equation}
g_{\lambda,\mu} \left( \left( \begin{array}{cr}
p_1 & -\bar{q}_1 \\ q_1 & r_1 \end{array} \right), \left(\begin{array}{cr}
p_2 & -\bar{q}_2 \\ q_2 & r_2 \end{array} \right) \right)= {\rm Re}
 (\lambda \bar{p}_1 p_2 + \bar{q}_1 q_2 + \mu \bar{r}_1 r_2)\, ,\label{eq5}
\end{equation}
for $\lambda, \mu$ real positive parameters, and where given any quaternion 
$x$, ${\rm Re}\, x$ denotes its real part. 

The norm defined by the metric $g_{\lambda,\mu}$ is given by 
$$
\left\| \left( \begin{array}{cr}
p & -\bar{q} \\ q & r \end{array} \right) \right\| ^2_{\lambda,\mu} =
\lambda|p|^2 + |q|^2 + \mu |r|^2\, .
$$
The bi-invariant metric on $Sp(2)$ is $g_{\frac{1}{2},\frac{1}{2}}$.

We summarize the properties that the family $g_{\lambda,\mu}$ has in 
the form of a proposition, whose verification should be a simple task
for the reader.

\begin{proposition} \label{pici}
Let $g_{\lambda,\mu}$ be the family of metrics on $Sp(2)$ defined by 
{\rm (\ref{eq5})} and let $g_{\lambda}$ be the induced metric on 
${\mathbb S}^7$. Then: 
\begin{enumerate}
\item \begin{equation}
\pi_{\circ}: (Sp(2),g_{\lambda,\mu}) \rightarrow ({\mathbb S}^7,g_{\lambda})
\label{sp2s}
\end{equation}
is a Riemannian submersion.
\item The fibers of this submersion are canonically placed totally geodesic 
$3$-spheres of constant sectional curvature $1/\mu$. 
\end{enumerate}
The family of metrics $g_{\lambda,\mu}$ is invariant under right 
multiplication by elements of the subgroup ${\mathbb S}^3 \times 
{\mathbb S}^3 \subset Sp(2)$, realized as the set of diagonal matrices.
\end{proposition}

Let $V_M^{\circ}$ and $H_{M}^{\circ}$ be the vertical and horizontal spaces 
of the Riemannian submersion $\pi_{\circ}$ in (\ref{sp2s}), 
respectively.
By (\ref{ac0}), if $L_{M}$ stands for the operator given by left 
translation by $M\in Sp(2)$, we see that 
\begin{equation}
V_{M}^{\circ}=L_{M*}\left\{ \left( \begin{array}{cc} 0 & 0 \\ 0 & p 
\end{array}\right): \; p \in {\mathbb H}, \; \bar{p}=-p\right\} \, .
\label{ve0}
\end{equation}
The space $H_{M}^{\circ}$ is the perpendicular of $V_{M}^{\circ}$ in the
$g_{\lambda,\mu}$-metric. This computation may be carried at the identity 
element
because the metric is left-invariant. Hence, by (\ref{eq5}), we see that
\begin{equation}
H_{M}^{\circ}=L_{M*}\left\{ \left( \begin{array}{cr} p & -\bar{q} \\ q & 0 
\end{array}\right): \; p,q \in {\mathbb H}, \; \bar{p}=-p\right\} 
\, .
\label{ho0}\end{equation}

We present expressions for $\nabla_{Z_0}Z_1$ and for the 
curvature $K_{\lambda, \mu}(Z_0,Z_1)$ of the planar section in
$T_M Sp(2)$ spanned by $Z_0$ and $Z_1$.
We consider the decompositions 
$$
Z_0=T+U\, , \quad Z_1=X+V
$$ 
of $Z_0$ and $Z_1$ in terms of their $\pi_{\circ}$ horizontal and vertical 
components, respectively. Since $g_{\lambda,\mu}$ is a left-invariant 
metric, we extend the vectors in a left-invariant manner, to reduce our
calculations to one nearby or at the identity, respectively. Hence, 
by (\ref{ve0}) and 
(\ref{ho0}), the vector $Z_0$ is of the form 
\begin{equation}
Z_0=T+U=L_{M*}\left(\begin{array}{rr}
p & -\bar{u} \\ u & 0 \end{array}\right)+
L_{M*}\left(\begin{array}{rr}
0 & 0 \\ 0 & r \end{array}\right)
 \label{ec0}
\end{equation}
for some quaternions $p$, $u$, $r$, with $p$ and 
$r$ purely imaginary, while the vector $Z_1$ is of the form  
\begin{equation}
Z_1=X+V = L_{M*}\left(\begin{array}{cr} 
q & -\bar{w} \\ w & 0 \end{array} \right) +
L_{M*} \left(\begin{array}{rr} 
0 & 0 \\ 0 & z \end{array} \right) \, , \label{ec1} 
\end{equation}
for some quaternions $q$, $w$ and $z$, with $q$ and $z$ purely imaginary
also. Though the formulas we provide hold in general, we could assume 
further that $Z_0$ and $Z_1$ are vectors
of norm one that are perpendicular to each other, facts encoded into the
relations
$$ 
\begin{array}{rcl}
\lambda|p|^2+|u|^2 +\mu |r|^2 & = & 1\, , \\ 
\lambda|q|^2+|w|^2+\mu |z|^2 & = & 1\, , 
\end{array}
$$ 
and
$$ 
{\rm Re}(\lambda\bar{p}q+\bar{u}w +\mu \bar{r}z)=0\, , 
$$  
respectively.

On any Lie group endowed with a left-invariant metric
$\< \, \cdot \, , \, \cdot \, \>$, we have that
\begin{equation}
\langle \nabla_U V, W\rangle =\frac{1}{2}(\langle [U,V],W\rangle -
\langle [V,W],U\rangle +\langle [W,U],V\rangle )\, , \label{cov}
\end{equation}
for $U$, $V$,  and $W$ arbitrary left-invariant fields. The 
 expression for the Levi-Civita connection of $g_{\lambda, \mu}$ follows
from it. Indeed, let us denote by 
${\mathcal V}$ and ${\mathcal H}$ the orthogonal projections onto the 
vertical and horizontal tangent spaces of $Sp(2)$, respectively. Further,
let us denote by ${\mathcal V}'$ and 
${\mathcal H}'$ the compositions of ${\mathcal H}$ with the 
projection onto the tangent to the Hopf fiber and the perpendicular 
component to it in ${\mathbb S}^7$, respectively. Then, 
\begin{equation}
\begin{array}{rcl}
\nabla_{Z_0} Z_1 \! \! \! \! \! & = & \! \! \! \! L_{M*}\! \!
\left( \begin{array}{cc} \frac{1}{2}(pq-qp-\bar{u}w+\bar{w}u) & \lambda
q\bar{u}+\mu \overline{zu}-(1-\lambda)p\bar{w}-(1-\mu) \overline{rw}  \\
\lambda uq-\mu zu -(1-\lambda)wp+(1-\mu)rw & 
\frac{1}{2}(w\bar{u}-u\bar{w}+rz-zr)
\end{array}\right) \vspace{1mm} \\
\! \! \! \! & = & \! \! \! \! 
L_{M*}\! \! \left( \frac{1}{2}[Z_0,Z_1] + (\lambda-\frac{1}{2})
\left( \begin{array}{cc}
 0 & q\bar{u}+p\bar{w} \\
uq+ wp & 0
\end{array}\right) 
-(\mu-\frac{1}{2})
\left( \begin{array}{cc}
 0 & \bar{u}z+\bar{w}r \\
zu+ rw & 0 
\end{array}\right) 
\right) 
\vspace{1mm} \\
\! \! \! \! & = & \! \! \! \!
L_{M*}\! \! \left( \frac{1}{2}[Z_0,Z_1]\! + \! (\lambda\!-\!\frac{1}{2})\left(
[{\mathcal H}'Z_0,{\mathcal V}'Z_1]\! +\! 
[{\mathcal H}'Z_1,{\mathcal V}'Z_0]\right)
\!+\!(\mu\!-\!\frac{1}{2})\left([{\mathcal H} Z_0,{\mathcal V}Z_1]\!+\!
[{\mathcal H}Z_1,{\mathcal V}Z_0]\right) 
\right)\, .
\end{array}
\label{cd}
\end{equation}

The expression for the sectional curvature of $g_{\lambda, \mu}$
is obtained using the fundamental formulas for Riemannian submersions 
of O'Neill \cite{on}. 
We state the result in the form of a theorem.

\begin{theorem}\label{Kur}
Consider $Sp(2)$ with the left invariant metric $g_{\lambda,\mu}$ 
defined by {\rm (\ref{eq5})}. 
Let $Z_0$ and $Z_1$ be left-invariant vectors fields
of the form {\rm (\ref{ec0})} and {\rm (\ref{ec1})}, respectively.
Then, the unnormalized 
curvature $K_{\lambda,\mu}(Z_0,Z_1)$ of the planar section of $T_{M}Sp(2)$ 
spanned by $Z_0$ and $Z_1$ is given by
\begin{equation}
\begin{array}{rcl}
K_{\lambda,\mu}(Z_0,Z_1) & = & 
\! (4-3(\lambda+
\mu))(|u|^2|w|^2-({\rm Re}(\bar{u}w))^2)+\lambda^2(|q|^2|u|^2+|p|^2|w|^2)
\\ & & +\lambda(|p|^2|q|^2-({\rm Re}(\bar{p}q))^2)-2{\rm Re}(\bar{u}w)\left(
\lambda{\rm Re}(\bar{p}q)+\mu{\rm Re}(\bar{r}z)\right)
 \\ & & +2\lambda(1-\lambda){\rm Re}(2\overline{wq}up-
\overline{uq}wp )+\mu^2(|u|^2|z|^2 
+|w|^2|r|^2) \\ 
& & +\mu(|r|^2|z|^2-({\rm Re}(\bar{r}z))^2)+
 2\lambda \mu {\rm Re}( (uq\bar{u}- up\bar{w})z +(wp\bar{w}-wq\bar{u})r)
\\ & & +2\mu(1-\mu){\rm Re}\left(2\overline{zw}ru-\overline{zu}
rw\right) \, .
\end{array}\label{rsec}
\end{equation}
When $\lambda=\mu=\frac{1}{2}$, this expression yields
$$
K_{\frac{1}{2},\frac{1}{2}}(Z_0,Z_1)=\frac{1}{4}\| 
[Z_0,Z_1]\|^2_{g_{\frac{1}{2},\frac{1}{2}}}\, ,
$$
and except for this one case, the
curvature operator ${\mathcal R}_{g_{\lambda,\mu}}$ admits negative 
eigenvalues. The scalar curvature of the metric is
$s_{g_{\lambda,\mu}}=2\left( \frac{3}{\lambda}+24-6(\lambda+\mu)+
\frac{3}{\mu}\right)$. For $(\lambda, \mu)\in (0,1/2]\times
(0,1/2]$, the sectional curvatures of $g_{\lambda, \mu}$ are nonnegative.
\qed
\end{theorem}

We may use the explicit expression for the sectional curvature given in Theorem
\ref{Kur} to prove that these curvatures are nonnegative for 
$(\lambda, \mu)\in (0,1/2]\times (0,1/2]$. This method is direct but cumbersome. 
Fortunately, earlier work of J. Cheeger \cite{jeff} may be applied to show this
 nonnegativity also in a more conceptual manner.
Cheeger's result pertains the study of Riemannian manifolds with isometric
group actions that shrink the metric in the direction of the orbits.

\begin{theorem}\label{t315}
Let $(\lambda, \mu)\in \mb{R}_{>0}^2$. Then all fibers of 
$\pi_{\circ}: (Sp(2),g_{\lambda,\mu}) \rightarrow 
({\mathbb S}^7,g_{\lambda})$ are critical points of the functional 
{\rm (\ref{sf})}, of critical value $0$. 
\end{theorem}

Let us observe that if we were to intertwine the 
parameters $\lambda$ and $\mu$ above in Proposition \ref{pici}, we 
would obtain analogous statements for
the fibration $\pi_{{\circ}_{1}}$ of (\ref{ac1}) instead. 
The volume of the fibers would then serve as a distinguishing criteria
to select the best among these two fibrations and their fibers.
We tie this issue with homology considerations next.

The integral homology of $Sp(2)$ is isomorphic to ${\mathbb Z}$ in
dimensions $0$, $3$, $7$, and $10$. Let $D$ be a positively oriented
generator of $H_3(Sp(2);{\mathbb Z})$. We consider the functional
(\ref{cpp}) in the case of this class. Being a class of 
dimension $3$, $D$ admits
realizations by submanifolds of $Sp(2)$ (cf. Remark \ref{thre}). In fact, 
all fibers of $\pi_{{\circ}_t}$, $t\in [0,1]$, realize $D$.

Using the techniques of Sullivan \cite{dennis2} and the calibrations of
Harvey \& Lawson \cite{hala}, Tasaki has shown \cite{tasa} that   
if $H$ is a closed Lie subgroup of a compact Lie group $G$ that represents
a nontrivial real homology class, then the group $G$ carries a left-invariant
metric such that the volume of $H$ in the induced metric is less than or equal
to the volume of any embedded representative of its homology class.
The following result follows by observing that in the case of 
$Sp(2)$ and the ${\mathbb S}^3$-fibers of $\pi_{\circ}$ above, 
any left-invariant metric on $Sp(2)$ 
for which the horizontal-vertical decomposition for $T_M Sp(2)$ is that
given by (\ref{ve0}) and (\ref{ho0}) works well in the alluded result.

\begin{theorem}\label{t316}
Let $D$ be the positively oriented generator of $H_3(Sp(2);{\mathbb Z})$. 
Assume that $\mu \leq \lambda $. Then in the induced metrics, all
fibers of  $\pi_{\circ}: (Sp(2),g_{\lambda,\mu}) \rightarrow 
({\mathbb S}^7,g_{\lambda})$
 are volume minimizer elements of
${\mathcal M}_D(Sp(2))$, of volume $2\pi^2 \mu^{\frac{3}{2}}$ each. 
\qed
\end{theorem}
\medskip

Hence, under the assumption that $\mu \leq \lambda$, the result above
makes of all the fibers of $\pi_{\circ}$ candidates for canonical 
representatives of $D$.

\begin{theorem}\label{t317}
Let $D$ be the positively oriented generator of $H_3(Sp(2);{\mathbb Z})$, and
let $(\lambda, \mu)\in (0,1/2]\times (0,1/2]$ with $\mu \leq \lambda$. Then 
the functional {\rm (\ref{cpp})} is bounded below 
by $0$, all fibers of $\pi_{\circ}: (Sp(2),g_{\lambda,\mu}) \rightarrow 
({\mathbb S}^7,g_{\lambda})$ realize its minimum value $0$ 
and are volume minimizer elements of ${\mathcal M}_D(Sp(2))$, and if
$(M,g) \in {\mathcal M}_D(Sp(2))$ is a volume minimizer $0$ of the
functional, then $M$ must be isometric to a fiber.
\end{theorem}

Thus, modulo isometries, the class 
$D$ has a unique canonical representative, a standard sphere of smallest
volume in $\mc{M}_D(Sp(2))$. 

{\it Proof}. By Proposition \ref{pici}, all fibers are totally geodesic
and therefore critical points of (\ref{sf}) and (\ref{cpp}), respectively,
of critical value $0$.
The fibers are all diffeomorphic to each other.
We show that an isometrically embedded manifold that represents
$D$ and is a volume minimizer minimum of (\ref{cpp}) 
must be isometric to a fiber of $\pi_{\circ}$.

Let $M\in \mc{M}_D( Sp(2))$ be totally geodesic and of minimal
volume. We ease into the general argument by first considering the case
where $\lambda=\mu=\frac{1}{2}$ so that the
background metric is the bi-invariant metric on $Sp(2)$. 
By (\ref{cd}), we see that the Lie bracket of left-invariant vector fields
tangent to $M$ is a left-invariant vector field tangent to $M$, and so the
restriction of the Lie bracket operation to $M$ gives $M$ the structure of 
a Lie subgroup of $Sp(2)$. Since the metric on $M$ has nonnegative sectional curvature, 
$M$ must be one of $\mb{S}^3$ or $\mb{S}\mb{O}(3)$ or $\mb{S}^1
\times \mb{S}^1\times \mb{S}^1$ as a closed Lie subgroup of $Sp(2)$, 
and the latter two are ruled out by the fact that $M$ is a homological
$3$-sphere (the $3$-torus can be ruled out also by using the
fact that the maximal torus in $Sp(2)$ is of dimension $2$; on the 
other hand, the orthogonal group can be ruled out by the fact that
the exponential in $Sp(2)$ of a $3$-dimensional Lie subalgebra isomorphic to 
$\mf{s}\mf{u}(2)$ is $\mb{S}^3$; the quotient of this universal cover
by $\mb{Z}/2$ yields $\mb{S}\mb{O}(3)\cong \mb{P}^3(\mb{R})$).
Thus, $M$ is in fact diffeomorphic to a $3$-sphere,
and since it has minimal volume among representatives of $D$, it must
be isometric to an $\mb{S}^3$ $\pi_{\circ}$-fiber.

In the general case when $Sp(2)$ is given the metric $g_{\lambda, \mu}$,
 $0<\mu \leq \lambda \leq 1/2$, the difficulty lies in showing that
a totally geodesic submanifold $M$ in $\mc{M}_D(Sp(2))$ of minimal
volume $2\pi \mu^{\frac{3}{2}}$ must be a Lie subgroup. We proceed as
follows.


For any $3$-manifold $(M,g_M)$ isometrically embedded into 
$(Sp(2),g_{\lambda,\mu})$, by the tubular neighborhood theorem, there exists
a neighborhood $\mc{O}$ of $M$ in $Sp(2)$ that is diffeomorphic to 
a neighborhood $N_{\mc{O}}$ of the zero section of the normal bundle of $M$ in 
$(Sp(2),g_{\lambda,\mu})$ via the exponential map. The bundle $N_{\mc{O}}$ can be
provided with two metrics, the first being just the pull-back of $g_{\lambda,\mu}$ on
$\mc{O}$ under the exponential map, and the second a metric compatible with the
normal connection $\nabla^{\nu}$ of the normal bundle $\nu(M)$ such that the projection
$\pi: N_{\mc{O}} \rightarrow M$ is a Riemannian submersion.
We let $A$ be the O'Neill tensor of this submersion. 

If $M$ is totally geodesic, the two metrics above on $N_{\mc{O}}$ are $C^1$-close and 
agree with each other to order $1$ on the zero section $M \hookrightarrow N_{\mc{O}}$
where they are just the background metric $g_{\lambda, \mu}$. Hence, on points of the 
totally geodesic zero section, the covariant derivatives that these metrics define 
coincide, and so do tensorial quantities defined in terms of the covariant derivatives and 
the metric tensor itself. 

For any Riemannian submersion, given any pair of horizontal vector
fields $X,Y$, we have that $A_X Y=\mc{V} (\nabla_X Y)=\mc{V}[X,Y]$. Here, $\nabla$ is the
covariant derivative of the Riemannian metric on the total space of the submersion, 
and $\mc{V}$ the projection onto the vertical space of the submersion. We apply this to the
Riemannian submersion above with O'Neill tensor $A$. Since $M$ is 
totally geodesic, if we consider left invariant vector fields $X$, $Y$ in $N_{\mc{O}}$ that 
are tangent to $M$ at one point $p\in M$, the tensorial quantity $(A_X Y)_p= 
\mc{V}(\nabla_X Y)_p$ may be computed by taking horizontal extensions $X$ and $Y$ 
of $X_p$ and $Y_p$, respectively,  and finding the projection onto the vertical component 
of $\nabla_X Y$.
But the extensions $X$ and $Y$ are then vector fields tangent to $M$ on points of $M$, and
along $M$ we have that $\nabla_X Y= \nabla^{g_M}_X Y$, which is tangent to $M$. Thus, 
$\mc{V}[X,Y]_p=0$ and so $[X,Y]_p$ is tangent to $M$. We conclude that the Lie bracket of 
left invariant vector fields that are tangent to $M\hookrightarrow N_{\mc{O}}$ at $p$ (and 
so horizontal at $p$) stay horizontal and are tangent to 
the zero section over it. This confers the tangent space of $M$ at any point with
the structure of a Lie subalgebra of the Lie algebra of $Sp(2)$, and up to diffeomorphism,
$M$ must be the image under the group exponential map of a $3$-dimensional Lie subalgebra
of $\mf{s}\mf{p}(2)$. Thus, $M$ is a closed Lie subgroup of $Sp(2)$. The remaining part of
the argument is the same as that for the bi-invariant metric given above.
\qed

The reasoning above can be extended to yield the following conclusion. 

\begin{corollary}
Consider the functional {\rm (\ref{cpp})} for the class $kD$, $k\in \mb{Z}$.
Its minimum is zero, and modulo isometries, the minimum of smallest volume is uniquely 
realized by a set of $k$ fibers if $k\geq 0$, or by $-k$ fibers each with 
the opposite orientation if $k<0$. Up to diffeomorphisms, 
all elements $[M]\in H_3(Sp(2);{\mathbb Z})$ admit a 
unique canonical representative, and the number of connected components of this
representative equals the divisibility $|[M]/D|\in \mb{Z}_{\geq 0}$ of the class $[M]$. 
\qed
\end{corollary}

\begin{remark}\label{bhk}
If we consider the $7$th homology group $H_7(Sp(2);{\mathbb Z})\cong {\mathbb Z}$ instead,
its generator $D$ cannot be realized by a submanifold \cite{bhk}. As ${\rm dim}\, Sp(2)=10$, 
this result sharpens  
Thom's original example of a nonrealizable $7$-class in a manifold
of dimension $14$ (see Remark \ref{thre}). 
Thus, for this homology class $D$ the space ${\mathcal M}_D(Sp(2))$ is empty.
Proceeding along similar lines, we may study a problem analogous to 
the one above by enlarging ${\mathcal M}_D(Sp(2))$ to include now all isometric 
immersed representatives of $D$, which for this $7$-class would be 
a nonemtpy space. We could attempt then to decide if the infimum of
the functional (\ref{cpp}) is realized over this extended domain, and perhaps tie up
the geometric properties of the minimizers with the nonrepresentability of the class.
\qed
\end{remark}

\section{Complex submanifolds of complex projective space}
\label{s4}
\setcounter{theorem}{0}
In this section, the ambient space $\tm$ is assumed to be ${\mathbb P}^n(\mb{C})$ 
provided with the Fubini-Study metric. We show 
that any of its complex submanifolds is a critical point
of (\ref{sf}). 

\subsection{The Fubini-Study metric}
We realize ${\mathbb P}^n({\mb{C}})$ as the quotient 
${\mathbb S}^{2n+1}/{\mathbb S}^1$, where the 
action of ${\mathbb S}^1$ on ${\mathbb S}^{2n+1}\subset {\mathbb C}^{n+1}$
is by complex scalar multiplication. The Fubini-Study metric $\tg$ on it
may be defined by the K\"ahler form of the metric, 
which on an affine chart $U$ with holomorphic coordinates 
$z=(z^1, \ldots ,z^n)$ is given by 
$$
\omega_U = \frac{i}{2} \ddb \log{(1\! +\!z^1 \bar{z}^1\! + \! \cdots
\! + \!  z^n
\bar{z}^n)}
=\frac{i}{2} \frac{\delta_{\a \b}(1\!+ \!\sum z^\gamma \bar{z}^{\gamma} )-
z^{\beta} \bar{z}^{\a}}{(1+\sum z^\gamma \bar{z}^{\gamma} )^2}dz^{\a}
\wedge d\bar{z}^{\b} \, .
$$
These locally defined forms give rise to
a globally defined form $\o$ that is $J$-invariant,
closed, and positive. 

The Ricci form and scalar curvature of the Fubini-Study metric are given
by $\rho_{\tg} = 2(n+1)\omega$ and $s_{\tg}=4n(n+1)$, 
respectively. 
The first Chern class of complex projective space is given by the positive
class $c_1({\mathbb P}^n({\mb{C}}))=((n+1)/\pi)[ \o]$.

Let $\bar{g}$ be the standard metric on ${\mathbb S}^{2n+1}$. With the
Fubini-Study metric $g$ as defined above on the base, we obtain
a Riemannian submersion ${\mathbb S}^1 
\hookrightarrow {\mathbb S}^{2n+1}\rightarrow {\mathbb P}^n({\mb{C}})$. 
Let $\{u,v\}$ be an orthonormal basis for a planar section 
at a point $p$ of $\mb{P}^n({\mb{C}})$, and 
let $\{\bar{u}, \bar{v}\}$ denote the horizontal lift of this basis to 
a point $\bar{p}$ of the fiber over $p$. Then we 
have \cite{on} that
\begin{equation}
K_g(u,v)=1+3|\bar{g}(\bar{u},J\bar{v})|^2\, , \label{sect}
\end{equation}
and so, the sectional curvature of the Fubini-Study metric 
ranges in the interval $[1,4]$, with the maximum attained by holomorphic 
or antiholomorphic sections, that is to say, sections spanned by $u,v$ where $v=\pm Ju$, while the
minimum is attained by sections spanned by $u,v$ where $Jv$ is orthogonal
to $u$ also. In fact, if $R^g$ is the $(0,3)$-curvature tensor of
$g$, we have that 
\begin{equation}
R^g(v,u)u=\cos{(t)}v_0 + 4\sin{(t)}Ju \, , \label{rc}
\end{equation}
where 
$$
v=\cos{(t)}v_0 + \sin{(t)}Ju 
$$
is the decomposition of $v$ with $v_0$ orthogonal to both, $u$ and $Ju$.


\subsection{Complex submanifolds of complex projective space}

\begin{proposition}\label{p11}
For any complex Riemannian manifold $(M,g)$ isometrically immersed in 
$\mb{P}^n({\mb{C}})$, the function $\| \a \|^2$ is an intrinsic
invariant, that is to say, it does not depend upon the immersion.
\end{proposition}

{\it Proof}. Let $m$ be the complex dimension of $M$. We use (\ref{sc}) and 
(\ref{sect}) to conclude that
\begin{equation}
s_g= 4n(n+1)-2(n-m)(4+2(n+m-1))- \| \a \|^2 =4m(m+1)-\| \a \|^2 \, . \label{ex}
\end{equation}
Since the scalar curvature is intrinsic, the result follows.
\qed

\begin{remark}
If the background manifold were to be the space
form $S^{\tn}_c$, or any of its quotients by a discrete subgroup, then for any 
isometrically immersed submanifold $(M,g)$ the quantity 
$\| \a\|^2 -\|H\|^2$ would be intrinsic also. This would follow by (\ref{sc}) 
using the same argument above, all the sectional curvatures now being the constant
$c$.
\qed
\end{remark}

\begin{theorem}\label{th7}
Any complex submanifold of $\mb{P}^n({\mb{C}})$ is a critical point 
of {\rm (\ref{sf})}. 
\end{theorem}

{\it Proof}. We use the critical point equations as described in 
Theorem \ref{th2}. By the minimality of every complex submanifold of a 
K\"ahler manifold, our task reduces to verifying that for any $m$,
$1\leq m \leq q$, we have that
\begin{equation} 
\<R^{\tg}(\a (e_i,e_j),e_j)e_i,\nu_m \> +   
\<\n_{e_i}^{\tg}(R^{\tg}(e_j,e_i)e_j)^{\nu}, \nu_m \> 
+{\rm trace}\, A_{\nu_m}A_{\nu_k}^2= 0 \, ,
\label{eql}
\end{equation}
which we do by verifying that each of the three summands on the left is
identically zero.

We decompose $\a(e_i,e_j)$ as $\a(e_i,e_j)=\< A_{\nu_l}e_i,e_j\>\nu_l$.
Let us consider a fixed index $l$ in the range $1\leq l\leq q$. We then
choose a basis $\{ e_i\}$ that diagonalizes $A_{\nu_l}$. Since
$\nu_l$ is 
orthogonal to both, $e_i$ and $Je_i$, by (\ref{rc}) 
we obtain that
$$
\< A_{\nu_l}e_i,e_j\>\<R^{\tg}(\nu_l,e_j)e_i,\nu_m \>=
\< A_{\nu_l}e_i,e_i\>\< \nu_l , \nu_m\> =0 \, ,
$$
where in deriving the last equality in the case when $l=m$,
we use the property that complex submanifolds of a K\"ahler manifold are 
austere. 
Thus, 
$$
\<R^{\tg}(\a (e_i,e_j),e_j)e_i,\nu_m\>= \< A_{\nu_l}e_i,e_j\>
\<R^{\tg}(\nu_l,e_j)e_i,\nu_m \>=0\, .
$$

In dealing with the middle summand, it is convenient to 
assume that the normal fields $\{ \nu_1, \ldots, \nu_q\}$ has been
chosen so that $(\n^{\tg}\nu_i)_p^\nu=0$. We then use (\ref{rc}) once
again to conclude that
$$
\begin{array}{rcl}
\<\n_{e_i}^{\tg}(R^{\tg}(e_j,e_i)e_j)^{\nu}, \nu_m \> & = & e_i\<
R^{\tg}(e_j,e_i)e_j)^{\nu}, \nu_m\> \\ & = & 
e_i\< R^{\tg}(e_j,e_i)e_j),\nu_m\> \\ & = & 0\, ,
\end{array}
$$
because, for each $i\neq j$, the decomposition of $e_i$ as 
$\cos{(t)}(e_i)_0 + 4\sin{(t)}Je_j$, with $(e_i)_0$ orthogonal to both,
$e_j$ and $Je_j$, produces a tangent vector to the complex submanifold.

Since for any normal field $\nu$ we have that 
$A_{\nu}J=-JA_{\nu}$, any of the cubic traces will have to be zero also. 
For in a
neighborhood of a fixed but arbitrary point $p\in M$, we
may choose an oriented local normal frame for $TM$ of the form
$\{ v_1, Jv_1, \ldots , v_d, Jv_d\}$, and obtain that
$$
\begin{array}{rcl}
{\rm trace}\, A_{\nu_m}A_{\nu_k}^2 & = & 
(\<A_{\nu_m}A_{\nu_k}^2 v_i,v_i\>+
\<A_{\nu_m}A_{\nu_k}^2 Jv_i,Jv_i\>) \\ & = & 
(\<A_{\nu_m}A_{\nu_k}^2 v_i,v_i\>
-\<J A_{\nu_m}A_{\nu_k}^2 v_i,Jv_i\>)
=0 \, .
\end{array}
$$
This completes the proof.
\qed
\medskip
 
\begin{remark}
The argument used above to show the vanishing of 
${\rm trace}\, A_{\nu_m}A_{\nu_k}^2$ holds for complex submanifolds of
an arbitrary background K\"ahler manifold. It is in the vanishing of the
other two terms that we used properties specific to the Fubini-Study metric on
$\mb{P}^n({\mb{C}})$.
It would be of interest to see the extent to which the stated result holds 
over arbitrary K\"ahler manifolds. A simple case to consider is that of the
product $\mb{P}^{n_1}({\mb{C}})\times \cdots \times \mb{P}^{n_k}({\mb{C}})$,
where the argument above can be adapted to prove that
\smallskip

\noindent {\bf Theorem $\ref{th7}'$.}
Any complex submanifold of the product $\mb{P}^{n_1}({\mb{C}})\times \cdots \times 
\mb{P}^{n_k}({\mb{C}})$
is a critical point of {\rm (\ref{sf})}. 
\qed

Using Corollary \ref{c313}, we derive a more general result:

\begin{theorem}
Let $(\tm ,J,\tg)$ be a K\"ahler manifold, and $M$ be a complex submanifold 
provided with the induced K\"ahler metric $g$. If the extrinsic scalar curvature of $M$ is
nowhere zero, then $M$ is a critical point of {\rm (\ref{cpp})} within its 
homology class. 
\end{theorem}

{\it Proof}. The total scalar curvature of $(M,g)$ is a topological invariant, and $(M,g)$ is
a critical point of (\ref{esc}) and (\ref{mc}). 
\qed    
\end{remark}
\medskip

\begin{remark}
There are critical points of (\ref{sf}) in 
$\mb{P}^n({\mb{C}})$ that are not complex submanifolds. For instance, the 
projection of the minimal torus
$\mb{S}^1(1/\sqrt{3})
\times  \mb{S}^1(1/\sqrt{3})\times \mb{S}^1(1/\sqrt{3})\subset \mb{S}^{5}$
onto $\mb{P}^2({\mb{C}})$ is a totally real minimal flat torus $T^2$ for which
$\| \alpha\|^2 =2$.
In affine holomorphic coordinates
$(z^1,z^2)$, $T^2=\{(z^1,z^2): \; |z^1|=|z^2|=1\}$. 
We may see that $T^2$ is a critical point of (\ref{sf}) in $\mb{P}^2(\mb{C})$
through a direct argument, showing that the summands in (\ref{eql}) are all zero. 
Clearly, this torus is homologically trivial in $\mb{P}^2({\mb{C}})$.  
\qed
\end{remark}
\medskip

Suppose that in the general setting we take the background manifold to be 
K\"ahler, say $(\tm,\tilde{J}, \tg)$. Let $D\in H_n(\tm ;{\mathbb Z})$. 
Since complex submanifolds minimize the volume within their homology
class, then we know that if 
the minimum of (\ref{cpp}) were to be realized by a complex manifold, 
such a critical point would be a critical point of
(\ref{sf}), (\ref{mc}), and of (\ref{ssf}), separately. In the case of
$\mb{P}^n({\mb{C}})$, 
by Theorem \ref{th7} we know that smooth algebraic subvarieties that 
represent elements of $H_{*}(\mb{P}^n({\mb{C}});{\mathbb Z})$ are
critical points of (\ref{sf}), 
and consequently, they are critical points of (\ref{ssf}) and (\ref{cpp})
also. 
We shall analyze later on if algebraic curves are, in fact, 
canonical representatives of their homology classes in 
$\mb{P}^2({\mb{C}})$. For now, we have the following:

\begin{theorem}\label{t44}
All complex submanifolds in 
${\mathcal M}_D(\mb{P}^n({\mb{C}}))$ are critical points of 
{\rm (\ref{sf})}, {\rm (\ref{ssf})} and {\rm (\ref{cpp})}, and they all 
have the same critical value.
\end{theorem}

{\it Proof}. By Theorem \ref{th7}, complex submanifolds are critical points
of these functionals. And for any two complex $m$-submanifolds in 
${\mathcal M}_D(\mb{P}^n({\mb{C}}))$, their volumes $[\omega]^m/m!$ and cup product
$c_1 \cup [\omega]^{m-1}$ coincide ($c_1$ being the first Chern class). 
Since $4\pi/(n-1)!$ times the latter is the total scalar curvature, 
the statement about the critical value follows by integration of the 
expression (\ref{ex}) in Proposition \ref{p11}.
\qed

We denote by $\mb{P}^n(\mb{R})$ the real projective $n$-space endowed with
its natural Riemannian metric. In general, if
$V$ is a real nonsingular variety, we denote by $V({\mathbb R})$ and 
$V({\mathbb C})$ its set of real and complex points, respectively. We
identify $V({\mathbb R})$ with the set of fixed points of complex 
conjugation in $V({\mathbb C})$. In this manner, $\mb{P}^n(\mb{R})$
is regarded as a real $n$-dimensional submanifold
of $\mb{P}^n(\mb{C})$. If $i: 
\mb{P}^n(\mb{R}) \hookrightarrow \mb{P}^n(\mb{C})$ is the inclusion map, then
$V(\mb{R}) \subset \mb{P}^n(\mb{R})  \hookrightarrow V(\mb{C}) \subset 
\mb{P}^n(\mb{C})$. 

\begin{corollary}
Any $k$-dimensional real subvariety $V=V^k$ of 
$\mb{P}^n(\mb{R})$ is a critical point 
of {\rm (\ref{sf})}. If we assume further that both $V(\mb{R})$ and
$\mb{P}^n(\mb{R})$ are oriented, with the orientations compatible with each 
other, and if we let $D=[V(\mb{R})]$ in 
$H_k(\mb{P}^n(\mb{R}); \mb{Z})$, then
all real varieties in ${\mathcal M}_D(\mb{P}^n(\mb{R}))$
are critical points of {\rm (\ref{ssf})} and 
{\rm (\ref{cpp})}, and they all have the same critical value.
\end{corollary}

{\it Proof}. Real projective space is a totally geodesic submanifold of 
complex projective space. 
\qed


\subsection{Algebraic curves in ${\mathbb C}{\mathbb P}^2$} 
A curve $S_d$ of degree $d$ in $\mb{P}^2(\mb{C})$ is, up to
diffeomorphisms, the Riemann surface given by the zeroes of the
polynomial $p_d(Z_0, Z_1, Z_2)=
Z_0^d + Z_1^d + Z_2^d$:
$$
S_{d}=\{ Z=[Z_0:Z_1:Z_2] \in {\mathbb C}{\mathbb P}^2: \;
p_d(Z)=0\}\, .
$$
A complex projective line $H$ in $\mb{P}^2(\mb{C})$ defines a
generator $h=[H]$ of the second integral homology group 
$H_2(\mb{P}^2(\mb{C});
{\mathbb Z})$, and $S_d$ represents the integer homology element
$dh$. Since the total Chern class of $\mb{P}^2(\mb{C})$ is given by
$c(\mb{P}^2(\mb{C}))=(1+\tilde{h})^3$, where $\tilde{h}=[\omega]/\pi$ 
is the generator of $H^2(\mb{P}^2(\mb{C}); {\mathbb Z})$ dual to 
the homology element $h$ above, if $i: S_d \hookrightarrow 
\mb{P}^2(\mb{C})$ is the inclusion map and
$x=i^{*}\tilde{h}$, we have that $c_1(S_d)=(3-d)x$. 

\begin{proposition}\label{p5}
The $L^2$-norm of the second fundamental form of a curve of degree $d$ in
$\mb{P}^2(\mb{C})$ is given by
$$
\Pi(S_d)=4\pi d(d-1)\, .
$$
\end{proposition}

{\it Proof}. In this dimension, the Fubini-Study metric has scalar curvature 
$24$. Since the total scalar curvature of a K\"ahler metric with K\"ahler 
class $\Omega$ is $4\pi c_1 \cdot \Omega^{n-1}/(n-1)!$, 
by (\ref{sc}) and (\ref{sect}) we see that
$$
4\pi c_1(S_d)\cdot [S_d]= 4\pi (3-d)d = 8\pi d -\int \| \a \|^2 d\mu \, ,
$$
and the result follows. 
\qed
\medskip

We may apply Gauss-Bonnet to rewrite $c_1(S_d) \cdot [S_d]$ in terms of
the Euler characteristic of $S_d$. Using the value of
$\Pi(S_d)$ computed above, (\ref{sc}) then yields that the 
genus $g_{S_d}$ of $S_d$ is given by
$$
g_{S_d}=\frac{(d-1)(d-2)}{2} \, .
$$
This genus of $S_d$ can be had via an argument independent of the one above,
using either the adjunction or Hurwitz's formula also.

We consider the integrand $\sum_{1\leq i,j\leq n}K_{\tg}(e_i,e_j)$ of the 
functional $\Theta$ in (\ref{esc}) when $M$ is a (real) surface $S$ in 
${\mathbb P}^2$ representing a homology class $D$.
We identify $S\in {\mathcal M}_D(\mb{P}^2(\mb{C}))$ with its
image under the embedding,  and let
$\{ e_1, e_2\}$ and $\{ \nu_1, \nu_2\}$ be orthonormal frames for its tangent 
and normal bundles, respectively. Thus, 
$\{ e_1, e_2, \nu_1, \nu_2\}$ is an orthonormal frame for the tangent 
bundle of the projective plane on points of $S$. The function 
$$ 
2K_{\tg}(e_1,e_2)\, .
$$ 
on $S$ has a maximum value $8$ at a point $p\in S$ if the frame   
$\{ e_1, e_2, \nu_1, \nu_2\}$ at this point is such that $Je_1=\pm e_2$ and 
$J\nu_1 = \pm \nu_2$. If the frame is positively oriented, this condition is 
equivalent to saying that the frame is holomorphic. On the other hand, the 
minimum of the function is
$2$, and this is achieved at points where the frame $\{ e_1, e_2, \nu_1, \nu_2\}$ is
purely real. Indeed, 
by (\ref{sect}), we have that
$$
2K(e_1,e_2)= 2+6|\bar{g}(\bar{e}_1,J\bar{e}_2)|^2\, ,
$$
where $\bar{g}$ is the standard metric on $\mb{S}^5$, and $\bar{e}_i$ is the
horizontal lift of $e_i$ in the Riemannian submersion 
${\mathbb S}^1 \hookrightarrow {\mathbb S}^{5}\rightarrow \mb{P}^2(\mb{C})$ 
to a point on the fiber over $p$. The statements made above follow readily from this.

We have the following:

\begin{theorem}
If $d\neq 3$, there are no Lagrangean elements in $\mc{M}_D(\mb{P}^2(\mb{C}))$.
\end{theorem}

{\it Proof}. Let $c_1$ be the pullback to $M$ of the first Chern class of $\mb{P}^2
(\mb{C})$.
Then (cf. with the proof of Proposition \ref{p5}) we have that 
$c_1 \cdot [M]=(3-d)d$. 
\qed 


\begin{theorem}
If $f_t : S_d \hookrightarrow \mb{P}^2(\mb{C})$ is a small deformation of 
an algebraic curve
$S_d $ in $\mb{P}^2(\mb{C})$ of degree $d$. Then  
$$
\Theta(f_f(S_d))=2\int_{f_t(S_d)}K(e_1,e_2) d\mu \geq 8\pi d \, ,
$$
$$
\mc{S}(f_t(S_d))=\int_{f_t(S_d)} (\| \alpha \|^2 - \| H\|^2) d\mu \geq 4\pi d(d-1) \, ,
$$
and 
$$
\Theta(f_f(S_d))- \mc{S}(f_t(S_d)) = 4\pi d(3-d) \, .
$$
\end{theorem}

{\it Proof}. The extrinsic scalar curvature of any algebraic variety $M$ of dimension $m$ in 
$\mb{P}^n(\mb{C})$ is given by
$$
\Theta (M)=4m(m+1)\int_M d\mu \, ,
$$
so its second variation under small deformations of the variety
is $4m(m+1)$-times the second variation of the volume form. 
By \cite[Proposition 3.2.2, Theorem 3.5.1]{si},
the latter is given by a nonnegative bilinear form on the 
normal bundle of $M$ whose kernel is the space of globally defined holomorphic cross sections 
of the normal bundle. Thus, under small variations of $M$, $\Theta(M)$ is either deformed into
a variety isotopic to $M$ with the same extrinsic scalar curvature or into a submanifold 
homologous to $M$ of larger extrinsic scalar curvature.    
The result now follows by the Gauss-Bonnet theorem, since (\ref{sc}) implies that
$$
\int_{f_t(S_d)} s_g d\mu_g = \Theta (f_t(S_d))-\mc{S}(f_t(S_d))\, ,
$$
and the Euler characteristic of any embedded surface stays constant under small 
perturbations of the surface. 
\qed

We shall improve the small deformation result above, and show that it 
holds for deformations in $\mc{M}_{[S_d]}$. 

\section{Oriented surfaces in $\mb{S}^2 \times \mb{S}^2$ and
$\mb{P}^2(\mb{C})$}
\label{s5}
\setcounter{theorem}{0}

\subsection{The diagonal in the product $\mb{S}^2 \times \mb{S}^2$}
\label{s51}
Let us examine the case of 
$\tm=\mb{S}^2\times \mb{S}^2$ provided with the product of the standard metrics and
 complex structures on the factors. The nontrivial 
integer homology groups are $\mb{Z}$, $\mb{Z}\oplus \mb{Z}$ and $\mb{Z}$ 
in dimensions $0$, $2$ and $4$, respectively. Let 
$A=[\mb{S}^2\times \{ q\}]$ and 
$B=[\{ p\}\times \mb{S}^2]$. Then $\{A,B\}$ 
generates $H_2(\tm;\mb{Z})$. 
We consider the problem of finding a canonical representative of
$D=A+B$. 

The union of the submanifolds $\mb{S}^2\times \{ q\}$ and $\{ p\}\times \mb{S}^2$ 
represents the class $D$. These submanifolds intersect at the point $p\times q$, 
and the intersection is modeled in complex coordinates $\{ z^1, z^2\}$ by the set 
$L=\{ (z^1,z^2)\in \mb{C}^2: \; z^1 z^2=0\, , \; | z^1 |^2 +| z^2 |^2 \leq 1\}$. 
The singularity of this representative can be resolved
by removing the pair $(B,L)$ in $\mb{C}^2$, where $B$ is the $z=(z^1,z^2)$-ball of radius
$1$ centered at the origin, and replacing it by 
the pair $(B,L')$, where $L'$ is a perturbation of 
$\tilde{L}=\{ (z^1,z^2): \; z^1 z^2= \varepsilon \, , 
\; | z^1 |^2 + |z^2 |^2 \leq 1\, , \; 0 < \varepsilon \ll 1 \}$ such that
$\partial L=\partial L'\subset \partial B$ (see \cite[p. 38-39]{gost}). This
procedure produces a smooth surface $M$ in $\mb{S}^2\times \mb{S}^2$ 
that represents $D$, $D=[M]$, because the sets $L$ and $L'$ are 
homologous to each other in $(B,\partial B)$. The surface $M$ is not 
a complex curve.  

The surface $M$ is homologous to the diagonal
$\mc{D}=\{ (p,p): \; p\in \mb{S}^2\}$, surface that represents 
the class $\mc{D}$ also. The diagonal $\mc{D}$ is a critical point of
(\ref{sf}), and being a K\"ahler submanifold of
the ambient K\"ahler manifold, it is a critical point of (\ref{mc}) also
of minimal volume among all the elements in $\mc{M}_D(\mb{S}^2 \times \mb{S}^2)$.

\begin{theorem}\label{th61} 
Let the product $\mb{S}^2 \times \mb{S}^2$ be provided with  
its standard metric and complex structure, and let us consider the
homology class $D=[\mb{S}^2\times \{ q\}]+[\{ p\}\times \mb{S}^2] \in 
H_2(\mb{S}^2 \times \mb{S}^2; \mb{Z})$. Let 
$\mc{D}$ be the diagonal in $\mb{S}^2 \times \mb{S}^2$.
Then $\mc{D}$ is a critical point of {\rm (\ref{sf})}, {\rm (\ref{mc})},  
{\rm (\ref{ssf})}, and {\rm (\ref{cpp})},
respectively, in each case, of critical value zero, and modulo isometries, 
$\mc{D}$ is the canonical representative of $D$. 
\end{theorem}

{\it Proof}. By Theorem $\ref{th7}'$, $\mc{D}$ is a critical point of
(\ref{sf}). In fact, $\mc{D}$ is a totally geodesic K\"ahler submanifold
of $\mb{S}^2 \times \mb{S}^2$ of scalar curvature $1$ that is a 
volume minimizer element of 
$\mc{M}_D(\mb{S}^2\times \mb{S}^2)$, of volume $8\pi$. 
We now prove that on $\mc{M}_D(\mb{S}^2\times \mb{S}^2)$, 
the functional (\ref{cpp}) is optimally bounded below by zero, and that
a minimum of smallest volume is a totally geodesic sphere of constant 
scalar curvature $1$.

Let $M$ be an oriented connected Riemannian surface $(M,g)$ isometrically 
embedded into $(\mb{S}^2 \times \mb{S}^2,\tg)$ that minimizes 
(\ref{cpp}) and have the smallest volume among minimizers. Then $M$ is a totally 
geodesic surface in $(\mb{S}^2 \times \mb{S}^2, \tilde{g})$ that is connected
 and has volume $8\pi$.    

Let $X,Y$ be an orthonormal 
basis of the tangent space $T_pM$ at $p\in M$. 
We have that $K_{\tg}(X,Y)=|X_1 \wedge Y_1|^2
+|X_2 \wedge Y_2|^2$, where $X=X_1+X_2$ and $Y=Y_1+Y_2$ are the decompositions 
of $X$ and $Y$ into components tangential to the factors, respectively. 

By (\ref{sc}) and the Gauss-Bonnet
theorem, we obtain that
$$ 
0=\int_M \| \alpha \|^2 d\mu = \int_M (\| \alpha\|^2-\| H\|^2) d\mu =
2\int_M K_{\tilde{g}}(X,Y) d\mu - 4\pi \chi(M) \, , 
$$
and since the sectional curvature of
the ambient space metric is nonnegative, $M$ can be either a torus or a 
sphere only. Let us view the product $\mb{S}^2\times \mb{S}^2$ as the total space of the
Riemannian submersion given by the projection onto one of the factors. 
If $M$ were a torus, then $K_{\tg}(X,Y)=0$ and at every $p \in M$ we could choose 
$\{ X,Y\}$ with $X$ tangent to the base and $Y$ tangent to the fiber. Using the 
exponential map,
we would then conclude that $M=\mb{S}^1 \times \mb{S}^1 \subset \mb{S}^2 \times
\mb{S}^2$, which contradicts the fact that $M$ represents a nontrivial class in
homology. So $M$ must be a topological sphere, $\chi(M)=2$, and       
$$
0=\int_M \| \alpha \|^2 d\mu =  \int_M (\| \alpha\|^2 -\| H\|^2) d\mu
= 2\int_M (|X_1 \wedge Y_1|^2 +|X_2 \wedge Y_2|^2) d\mu - 8\pi \, .
$$

Since horizontal geodesics stay horizontal, at no $p\in M$ it is possible to have
$T_p M$ horizontal or vertical. Otherwise, $M$ would be either a horizontal 
or a fiber sphere, respectively, neither one of which represents the homology
class $D$. Therefore, $M$ is the graph of a diffeomorphism $f: \mb{S}^2 \rightarrow
\mb{S}^2$, and by the transivity of the action of the group of isometries of    
the standard metric on $\mb{S}^2$, we conclude 
that $K_{\tilde{g}}(X,Y)=|X_1 \wedge Y_1|^2 +|X_2 \wedge Y_2|^2$ is a constant 
function on $M$, the constant $\frac{1}{2}$, that $f$ is an isometry of the
sphere with its standard metric, and that $M$ (the graph of $f$) is isometric to
$\mc{D}$ with its induced metric.  
\qed

The canonical representative $\mc{D}$ of $D$ above is a rigid complex 
submanifold of 
$\mb{S}^2 \times \mb{S}^2$ of zero genus in $\mc{M}_D(\mb{S}^2 \times \mb{S}^2)$. 
The immersed representative $\mb{S}^2\times \{ q\} \cup \{ p\} \times \mb{S}^2$ 
is an absolute minimizer of $\Pi$ also. Its local desingularizations at the 
intersection point produces a one parameter family of embedded spheres that 
represents $D$ and have a value of $\Pi$ close to but larger than $8\pi$, which
when subtracting the total scalar curvature of the representative yields a value
close to but larger than $0$.    
Thus, the functional $\mc{T}$ can be extended from $\mc{M}_D(\mb{S}^2 \times
\mb{S}^2)$ to the space of immersed representatives of $D$, and on this 
extended domain, there are singular current points of minimal mass that are 
absolute minimizers of the extension. The canonical representative $\mc{D}$ is 
a smooth current that lies in $\mc{M}(D)$, has $\Pi=0$, and so it minimizes the 
extension also, is a volume minimizer among representatives of $D$.
 
\begin{remark}
If $J'$ is any almost complex structure on $\mb{S}^2 \times \mb{S}^2$ tamed by
$\omega_1 + \omega_2$ (for instance, the $J$ in Theorem \ref{th61}), 
Theorem 2.4.C in \cite{gr} ensures the existence of a connected regular $J'$ 
curve into $\mb{S}^2\times \mb{S}^2$ that represents $A+B$ and has genus $1$, 
and the collection $M(J')_{1,1}$ of all such curves form a smooth manifold of 
dimension $6$. So for $J'=J$ these curves are ``symplectic canonical 
representatives'' of $A+B$, but none of them are metric 
canonical representative in the sense described by the Theorem above.    
\end{remark}

The class $2D=2(A+B)$ does not
admit an embedded sphere representative \cite{kemi}. This new class may be 
represented by the union of two horizontal and two fiber vertical
spheres, respectively, 
with this immersed representative an absolute minimizer of $\Pi$. 
Desingularizing at the four intersection points while preserving the homology 
class, we obtain an embedded noncomplex torus representative of $2D$ for which 
$\Pi$ is close to but larger than $16\pi$. Thus, when extended to the space of
immersed representatives, the value of $\Pi$ for the elements of the embedded 
family of representatives whose limit is the immersed one are larger than the 
value value zero of $\Pi$ for the limit.  The class $2D$ admits an 
$8$-parameter family of complex tori that represent it, each one of which is 
a volume minimizer within the class, of 
volume $4\pi$. We conjecture that these tori are the 
canonical representatives of $2D$, their second fundamental forms all having 
global $L^2$-norm equal to $16\pi$. This would show that we can find singular 
currents that are absolute minima of the extension of $\mc{T}$ to the space 
of immersed representatives of the class $2D$, which minimize the volume also. 
The gap between the functional values at this singular currents and the values 
at the smooth tori representatives would measure the minimum number of 
self-intersection points that immersed smooth representatives could have.   
The conjectured behaviour is reminescent of one parameter families of 
mappings that fail to be Fredholm just at a discrete set of parameters, where 
the mappings have finite dimensional kernel and cokernel but not closed 
range, and the index jumps \cite{simi}. 
The minimum value of
$\Pi$ on embedded manifolds reprentating an integral homology class jumps when 
the integral class changes.     


\subsection{Complex curves in $\mb{P}^2(\mb{C})$} 

In Theorem \ref{t44} we have shown that all complex submanifolds in  
${\mathcal M}_{[D]}(\mb{P}^n(\mb{C}))$ are critical points of the functional
$\mc{T}$ in (\ref{cpp})
of the same critical value, and by their minimality, 
also critical points of the domain restricted functional 
$\mc{S}\mid_{{\mathcal M}_{D}( \mb{P}^n(\mb{C}))}$ in (\ref{ssf}), all 
associated with the same critical value. It is then natural to ask if they 
are the minimum.  We study this optimization problem when $n=2$ and 
$D=[S_d]\in H_2(\mb{P}^2(\mb{C});\mb{Z})$ next. This will yield our first examples
of infimums of $\mc{T}$ over $\mc{M}_{[S_d]}(\mb{P}^2(\mb{C}))$ 
thar are no longer realized by totally geodesic
absolute minimizers of the functional. Thus, we are forced to study precompactness
properties of minimizing sequences.

Consider a minimizing sequence $\{M_n\}\subset \mc{M}_{[S_d]}(\mb{P}^2(\mb{C}))$:
$$
\mc{T}(M_n) \rightarrow  t=\inf_{M\in \mc{M}_{[S_d]}}\int_M  \| \alpha \|^2 
d\mu_M \, .
$$
We may assume that:
\begin{enumerate}
\item[M1.] \label{m1} All the $M_n$s are connected, and  $t$ is at most equal to
$\mc{T}(S_d)=4\pi d(d-1)$.
\item[M2.] \label{m2}
The $M_n$s can be taken to be of genus $g_{M_n}$ and areas $\mu_n$ satisfying the
bounds 
\begin{equation} \label{mg2}
g_{M_n} \leq \frac{(d-2)(d-1)}{2}+\frac{3}{4}d\, , 
\end{equation}
and 
\begin{equation} \label{ma2}
 \pi d \leq \mu_n \leq 2\pi(d(d-1)+2)\, .
\end{equation}
Indeed, by (\ref{sc}) and the bounds for the sectional curvature of the
the Fubini-Study metric in $\mb{P}^2(\mb{C})$,  we have that 
$$
\int_M \| \alpha \|^2 d\mu =
\int_M (2K_{g}(e_1,e_2)+\|H\|^2) d\mu
-4\pi \chi(M)\geq 2\mu(M) -4\pi\chi(M)\, ,
$$
which by M1 implies the upper bounds for $g_{M_n}$ and $\mu_n$, respectively.
On the other hand, $\mu(S_d)=\pi d$ is the absolute minimizer of
the area functional in $\mc{M}_{[S_d]}(\mb{P}^2(\mb{C}))$.
\item[M3.] By possibly passing to a subsequence, we may them   
assume that the minimizing sequence $\{M_n\}\subset \mc{M}_{[S_d]}(\mb{P}^2
(\mb{C}))$ is given by a sequence of isometric embeddings
\begin{equation} \label{ms}
f_n : (\Sigma, g_n) \rightarrow f_n(\Sigma)=M_n \hookrightarrow (\mb{P}^2(\mb{C}),g)
\end{equation}
of a fixed connected 
surface $\Sigma$, whose genus $g_{\Sigma}$ satisfies (\ref{mg2}), 
and where the areas $\mu_{g_n}(\Sigma)$ satisfy (\ref{ma2}). Notice that if 
$M=M_n=M_1 \# M_2$ where $M_2$ bounds some $3$-manifold in $\mb{P}^2(\mb{C})$, 
then $M_1$ must be homologous to $S_d$. We shall prove that $M_1$ is in effect
homeomorphic to $S_g$ and that $(M_n,g_n)$ converges to a complex curve of degree
$d$ endowed with the metric induced by the Fubini-Study metric on
$\mb{P^2}(\mb{C})$.
\item[M4.] We have the $L^1$ bound
\begin{equation}\label{pes2}
\left| \int_M s_g d\mu_g \right| \leq  \int_M | s_g | d\mu_g \leq 
\Theta(M) + \mc{T}(M) \, ,
\end{equation}
for the scalar curvature of all the $M_n$s in the minimizing sequence. This 
follows by (\ref{sc}), and the the Cauchy-Schwarz and Young inequalities, from
which we can derive the pointwise estimate 
\begin{equation}\label{pes}
| s_g | \leq 2K_{\tilde{g}}(e_1, e_2) + \| \alpha \|^2 
\end{equation}
for the scalar curvature of any $M\in \mc{M}_{[S_d]}(\mb{P}^2(\mb{C}))$. Then
(\ref{pes2}) follows by integration over $M$. 
\item[M5.] The minimizing sequence of isometric embeddings (\ref{ms}) can be 
taken to be given by conformal mappings with the scalar curvature of the metrics 
uniformly 
bounded in $L^2(\Sigma, d\mu_{g_n})$. 

{\it Proof}. We consider the conformal class defined by $g_n$ on $\Sigma$, and its
almost complex structure $J_n$. For dimensional reasons,  
the triple $(\Sigma, g_n, J_n)$ is K\"ahler. On such a manifold, 
there exists a conformal deformation $c_n: \Sigma \rightarrow
\Sigma$ to an extremal metric $c_n^* g_n$ of constant scalar curvature 
$\overline{s}_{g_n}=4\pi \chi(\Sigma)/\mu_{g_n}(\Sigma)$ \cite{ca}.
This extremal metric minimizes the $L^2$ norm of the scalar curvature of metrics 
on
the given conformal class. The result follows, as the areas $\mu_{g_n}(\Sigma)$ 
satisfy (\ref{ma2}). 
\qed 

Given a Riemannian surface $(M,g)$ isometrically embedded into 
$\mb{P}^2(\mb{C})$, we let 
$$
\mc{E}_g(M)=\int_M \| \alpha \|^4 d\mu_g 
$$
denote the ``energy'' of $M$. The statement above makes it natural to assume,
which we do, that the sequence $\{\mc{E}_{g_n}(M_n)\}$ is uniformly bounded.   
\end{enumerate}

Let us consider a minimizing sequence $\{(M_n,g_n)\}$ satisfying conditions 
M1-M5. In the conformal class of each element $(M_n,g_n)$ in this sequence,  
we shall consider a sequence of metrics $\{ g^c_m \}$ that
minimizes the energy $\mc{E}_{g}(M_n)$.
Notice that the conformal invariance of $\mc{T}(M_n)$ and $\mu_n$, and the 
Cauchy-Schwarz inequality imply that
$$
\mc{E}_g(M_n) \geq b_n:=\frac{\mc{T}^2(M_n)}{\mu_{g_n}}
$$
for any metric $g$ conformally equivalent to $g_n$, with equality if, and only 
if, $\| \alpha \|^2$ is constant. We call $b_n$ the energy
barrier of the conformal class of $g_n$. We study the
formation of singularities of $\{ g_m^c\}$ as $\mc{E}_{g_m^c}(M_n)\rightarrow
b_n$, and then study how these singularities evolve as the barrier $b_n$ varies 
with $g_n$ so that $\mc{T}(M_n) \rightarrow t$.
 These two steps are analogous to those in the study of
strongly extremal metrics \cite{sim,ss1,ss2},
 where we first study the existence of extremal metrics within a K\"ahler class, 
and then study how these metrics change as the value of the energy of one 
such varies as a function of the class towards the global minimum of the 
energy within the base of the K\"ahler cone. 
  
\subsubsection{Precompactness properties of the energy in a fixed conformal class}
Given a Riemannian surface $(M,g)$ isometrically embedded into 
$\mb{P}^2(\mb{C})$, we study the problem of minimizing $\mc{E}_g(M)$ within the
conformal class of $g$. In studying the formation of singularities of 
energy minimizing sequences, we follow closely the work of X.X. 
Chen \cite{xxc}, with the appropriate modifications required by our problem.

For any $(M,g)\in \mc{M}_{[S_d]}(\mb{P}^2(\mb{C}))$ that passes through a point 
$p$ in the projective plane, we define $D(p)$ to be
$$
D(p)=\{ q \in M: \; {\rm dist}\, (q,p)\leq \iota \} \, , 
$$
where $\iota=\frac{1}{2}{\rm inj\,  rad}_M$. 
If $r<\iota$, we denote by $D_r(p)$ the disk of radius $r$ centered at $p$.
We use conformal coordinates $z=x+iy$ on $D(p)$, and let $g_0=dx^2+dy^2$ be the
standard flat metric with Laplacian $\Delta_0$ and area form $d\mu_0$. Then the 
metric $g$ in $D(p)$ is of the form $g=e^{\varphi}g_0$ for some real valued 
function $\varphi$, thus identifying $g$ with $\varphi$, a fact that we use freely
below. We have the relation
\begin{equation} \label{csc}
s_g=e^{-\varphi } \Delta_0 \varphi \, ,
\end{equation}
where $s_g$ is the scalar curvature of $g$.

Given constants $C_1$, $C_2$, we define 
\begin{equation}
\mc{M}^{C_1,C_2}=
\{ \tg : \; \tg=e^\varphi g , \; 
\mu_{\tg}(M)\leq C_1\, ,  \mc{E}_{\tg}(M)\leq C_2 \} \, . 
\label{est}
\end{equation}
If $\tg \in  \mc{M}^{C_1,C_2}$ and $U$ is any
subdomain in $D(p)$, we define the local area, density and energy 
over $U$ by   
$$
\mu_{\varphi}(U)=\int_U e^{\varphi } d\mu_0 \, , \quad
t_{\varphi}(U)=\int_U \| \alpha \|^2  e^{\varphi } d\mu_0 \, , \quad 
e_{\varphi}(U)=\int_U \| \alpha \|^4  e^{\varphi } d\mu_0 \, , 
$$ 
and the local conformal invariant 
$$
\theta_{\varphi}(U)=\int_U 2K_{g}(e_1,e_2) e^{\varphi } d\mu_0 \, , 
$$ 
respectively. We have that 
$$
2\mu_{\varphi}(U) \leq \theta(U)\leq 8 \mu_{\varphi}(U) \, .
$$



Let $g_n$ be an energy minimizing sequence of metrics 
in the conformal class of $g$. 
We say that $p$ is a bubble point for $g_n$ if there exist constants $a$ 
and $e$ such that, for any $r$ with $0< r < \iota $, we 
have that  
$$
\liminf_n \mu_{\varphi_n}(D_r(p)) \geq a > 0 \, , \quad 
\liminf_n e_{\varphi_n}(D_r(p))\geq e \geq 0 \, , 
$$ 
respectively. The largest constants $a_p$ and $e_p$ for which this holds are the 
concentration weights of area and energy at the bubble point $p$.    

The local area and energy function at $p$ are defined by the expressions 
$$
A_p(r)= \limsup_n \mu_{\varphi_n}(D_r(p)) \, , \quad 
E_p(r)= \limsup_n e_{\varphi_n}(D_r(p)) \, , 
$$ 
and the local concentration of these quantities at $p$ are defined as the limits
$$
A_p= \lim_{r\rightarrow 0} A_p(r)\, , \quad 
E_p=\lim_{r\rightarrow 0}E_p(r) \, . 
$$ 
We say that $p$ is a basic bubble point if $A_p>0$, and $E_p \geq 0$. 
Passing to a subsequence if necessary, we see that every basic bubble point is a 
bubble point. 
 
The following isoperimetric inequality theorem plays for us the same 
essential role it has in the work of Chen \cite{xxc}. We state it here for 
completeness, and refer the reader to \cite{buza} for details and historical references
about it. This result will be used freely below.

\begin{theorem} \label{iso}
{\rm \cite{buza}}
Let $g$ be a metric in an Euclidean disk $D$ whose scalar curvature $s_g$ is in
$L^1(D,d\mu_g)$.
Then given any relatively compact disk $D'$ in $D$, we have that  
$$
\int_{D'} | s_g | d\mu_g \geq 4\pi 
-\frac{|\partial D'|_g^2}{\mu_g(D')} \, .
$$
\end{theorem}

We now have the following key result, the analogue of \cite[Lemma 2]{xxc}. 
It implies the existence of at most finitely many bubble points for sequences
in $\mc{M}^{C_1,C_2}$.  

\begin{lemma} \label{le54}
Let $\{g_n\}\in \mc{M}^{C_1,C_2}$ be an energy minimizing sequence in 
the conformal class of the metric $g$ on $M$. 
If $p$ is a bubble point for the $g_n$s, then 
$$
A_p(64A_p + E_p) \geq 8\pi^2 \, .
$$  
\end{lemma}

{\it Proof}. Given $\varepsilon >0$ sufficiently small, since $A_p(r)$ is a 
monotonically increasing function, we can choose a small enough radius $i_0$ such 
that  
$$
A_p \leq A_p(i_0)=\limsup_{n} \mu_{\varphi_n}(D_{i_0}(p)) < 
\left( 1+\frac{\varepsilon}{2} \right) A_p \, ,
$$
and then choose $n$ large enough so that
$$
\mu_{\varphi_n}(D_{i_0}(p))< (1+\varepsilon)A_p \, .
$$

For any $i_1$ such that $0<i_1 < i_0$, we consider Chen's  
waist concentration function $l_p(i_1, i_0)$:
$$
l_p(i_1,i_0)=\liminf_{n}\min_{i_1 \leq r \leq i_0} | \partial D_{r}(p)|_{g_n}
$$ 
By \cite[Lemma 1]{xxc}, we have that $\lim_{i \rightarrow 0}l_p(i,i_0)=0$. Thus,
we can choose $i_1 < i_0$ small enough such that $l_p(i_1,i_0)< \varepsilon$, an
integer $N=N(\varepsilon)$ such that 
$$
\min_{i_1 \leq r \leq i_0}| \partial D_r(p)|_{g_n} < 2\varepsilon \quad
\text{for all $n\geq N$,} 
$$
and $r_n \in [i_1, i_0]$ such that
$$
|\partial D_{r_n}(p)|_{g_n}< 3\varepsilon\, .
$$ 
It follows that 
$$
A_p \leq \mu_{\varphi_n}(D_{i_1}(p))\leq \mu_{\varphi_n}(D_{r_n}(p))\leq 
(1+\varepsilon) A_p \, ,
$$
and, by Theorem \ref{iso}, that
$$
\int_{D_{r_n}(p)} | s_{g_n} | d\mu_{g_n} \geq 4\pi 
-\frac{|\partial D_{r_n}(p)|_{g_n}^2}{\mu_{g_n}(D_{r_n}(p))} > 4\pi -
\frac{9\varepsilon^2}{A_p}>0\, . 
$$
By integration of the pointwise estimate (\ref{pes}), and the Cauchy-Schwarz
inequality, we conclude that
$$
\left(\int_{D_{r_n}(p)}|s_{g_n}| d\mu_{g_n} \right)^2 \leq  
(\theta_{\varphi_n}(D_{r_n})+t_{\varphi_n}({D_{r_n}}))^2 \leq 2( 
64\mu_{\varphi_n}^2({D_{r_n}}) +
\mu_{\varphi_n}({D_{r_n}})e_{\varphi_n}({D_{r_n}})) \, , 
$$  
and so we have that
$$
\frac{\left( 4\pi - \frac{9\varepsilon^2}{A_p}\right)^2}{(1+\varepsilon)A_p}
\leq \frac{\left(\int_{D_{r_n}(p)}|s_{g_n}| d\mu_{g_n} \right)^2 }{\mu_{\varphi_n}(
D_{r_n}(p))}\leq  2( 
64\mu_{\varphi_n}({D_{r_n}}) + e_{\varphi_n}({D_{r_n}})) 
\, .  
$$  
Since $r_n \leq i_0$, we conclude that
$$
\frac{\left( 4\pi - \frac{9\varepsilon^2}{A_p}\right)^2}{(1+\varepsilon)A_p}
\leq  2( 
64\mu_{\varphi_n}({D_{i_0}}) +
 e_{\varphi_n}({D_{i_0}})) \, , \quad n\geq N\, , 
$$  
which proves that
$$
2\liminf_n (64\mu_{\varphi_n}(D_{i_0}(p))+
e_{\varphi_n}(D_{i_0}(p))) >  
\frac{\left( 4\pi - \frac{9\varepsilon^2}{A_p}\right)^2}{(1+\varepsilon)A_p}\, .
$$
If we then let $\varepsilon \rightarrow 0$, we conclude that
$$
\liminf_n (64\mu_{\varphi_n}(D_{i_0}(p))+
e_{\varphi_n}(D_{i_0}(p))) \geq 
\frac{8\pi^2}{A_p}\, ,
$$
and so 
$$
64A_p(i_0)+E_p(i_0) \geq \frac{8\pi^2}{A_p} \, .  
$$  
The desired results follows by computing the limit as $i_0 \rightarrow 0$.
\qed

The following two lemmas of X. Chen \cite{xxc} hold in our context. We provide
a detailed proof for the second one of these only. Notice that   
the validity of the first requires only to have upper bounds control on the 
sequence of areas of the $M_n$s.

\begin{lemma} {\rm (\cite[Lemma 3]{xxc})} \label{l55}
Let $\{g_n\}\in \mc{M}^{C_1,C_2}$ be an energy minimizing sequence in 
the conformal class of the metric $g$ on $M$.  Then, 
$$
\liminf_{n} \min_{0\leq i_1 \leq r \leq i_0} \frac{1}{2\pi}\int_0^{2\pi} \varphi_n
(r\cos{\theta},r\sin{\theta}) d \theta 
$$
is bounded above for any $[i_1, i_0]$.  
\end{lemma}

\begin{lemma} {\rm (\cite[Lemma 5]{xxc})} \label{l56}
Let $\{g_n\}\in \mc{M}^{C_1,C_2}$ be an energy minimizing sequence in 
the conformal class of the metric $g$ on $M$.
If $p$ is not a basic bubble point for the $g_n$s {\rm (}that is to say, 
if $A_p=0${\rm )} and $i_0>0$ is sufficiently small,
there exists a constant $C>0$ and an integer $N$ such that
for $0<\alpha < 2$ and
$r\leq i_0$, we have that
$$
\frac{1}{r^{\alpha}}\mu_{\varphi_n}(D_r(p)) < C \, . 
$$
for all $n\geq N$.
\end{lemma}

{\it Proof}. We can choose 
$i_0$ sufficiently small such that $2A_p(i_0)(64 A_p(i_0)+E_p(i_0)) < 
(2\pi(2-\alpha))^2$.
Given this choice, let $C$ be a constant such that
$$
\frac{\mu_{\varphi_n}(D_{i_0}(p))}{i_0^{\alpha}} < C  \quad \text{for all $n$.}
$$
The desired statement holds for this value of $C$. 

Indeed, let us 
consider the function $f_n(r)=\mu_{\varphi_n}(D_{r}(p))-Cr^{\alpha}$. Assume that for 
all $n$ we have  
$i_n <i_0$ such that $f_n(i_n)>0$. Then there must exist $r_n\in (i_n,i_0)$ such that
$f_n(r_n)=0$ and $f'_n(r_n)<0$:
$$
\mu_{\varphi_n}(D_{r_n}(p))=Cr_n^{\alpha}\, , \quad
\int_0^{2\pi} e^{\varphi_n}r_n d\theta < C\alpha r_n^{\alpha-1}\, .
$$  
We then conclude that
$| \partial D_{r_n}(p)|_{g_n}^2 < 2\pi \alpha C r_n^{\alpha}$ because
$$
\left( \int_0^{2\pi} 
e^{\frac{1}{2}\varphi_n}r_n d\theta\right)^2 < 
\left( \int_0^{2\pi} 
e^{\varphi_n}r_n d\theta\right) \int_0^{2\pi} r_n d\theta 
< 2\pi C\alpha r_n^{\alpha}\, ,
$$
and by Theorem \ref{iso}, we see that
$$
\int_{D_{r_n}} |s_{g_n}|d\mu_{\varphi_n} > 4\pi - 
\frac{| \partial D_{r_n}(p)|_{g_n}^2}{\mu_{\varphi_n}(D_{r_n}(p))}
> 4\pi - 
\frac{2\pi C\alpha r_n^{\alpha}} 
{Cr_n^{\alpha}}=2\pi(2-\alpha) \, . 
$$ 
This contradicts the fact that since $r_n \leq i_0$, by (\ref{pes}), the 
square of the 
left side of this expression is 
bounded above by $2A_p(i_0)(64A_p(i_0)+E_p(i_0))$, and this quantity is 
bounded above 
by $(2\pi(2-\alpha))^2$. 
\qed

The following results will be used in combination with Lemma \ref{le54} to
show that sequences in $\mc{M}^{C_1, C_2}$ have at most finitely many bubble 
points.  
 
\begin{lemma}
Let $\{g_n\}\in \mc{M}^{C_1,C_2}$ be an energy minimizing sequence in
the conformal class of the metric $g$ on $M$.
If $p$ is not a basic bubble point for the $g_n$s, there exists a small
neighborhood $\mc{O}_p$ of $p$ and a positive constant $C$ such that
$$
\sup_{n} \sup_{q\in \mc{O}_p} \varphi_n (q) < C \, . 
$$
\end{lemma}

{\it Proof}.  If the the statement were false, we could modify the sequence 
of metrics slightly so that $\varphi_n(p)\rightarrow \infty$. We contradict this
fact. 

Let us choose a small enough disk $D_{i_0}(p)$ so that
Lemma \ref{l56} applies for some constant $C$ and integer $N$. 
 Consider the function $[i_1, i_0]\ni 
r\rightarrow \varphi_n (r\cos{\theta},r\sin{\theta})$, and for simplicity, set
$q_r(\theta)=(r\cos{\theta},r\sin{\theta})$.
Then we have that
$$
\begin{array}{rcl}
{\displaystyle 
\left| \int_{0}^{2\pi } (i_0 \partial_r \varphi_n(q_{i_0}(\theta)) d\theta
- \int_{0}^{2\pi } (i_1 \partial_r \varphi_n(q_{i_1}(\theta))d\theta \right|  
 }& = &  {\displaystyle \left|  \int_{i_1}^{i_0} \int_0^{2\pi}
(\partial_r( r \partial_r  \varphi_n (q_r(\theta)) d\theta dr  \right| }
\\ & = & {\displaystyle \left| \int_{i_1}^{i_0}\int_{0}^{2\pi}
 (r\partial_r^2 \varphi_n + \partial_r \varphi_n) d\theta dr  \right| } 
\\ & = & {\displaystyle \left| \int_{i_1}^{i_0}\int_{0}^{2\pi}
 r \Delta_0 \varphi_n  d\theta dr  \right| } 
\\ & \leq  & {\displaystyle \int_{i_1}^{i_0}\int_{0}^{2\pi}
 | s_{g_n} | d\mu_{g_n}  } \, , 
\\ & \leq  & 2^{\frac{1}{2}}(64\mu^{2}_{g_n}(A_{i_0,i_1}^p)
+ \mu_{g_n}(A_{i_0,i_1}^p)E_{g_n}(A_{i_0,i_1}^p))^\frac{1}{2}\, , 
\end{array}
$$
where $A_{i_0,i_1}^p$ is the annular region $D_{i_0}(p)\setminus D_{i_1}(p)$. In 
obtaining the last two inequalities we have used (\ref{csc}) and  (\ref{pes}),
respectively. Passing to the limit as $i_1\rightarrow 0$, and setting $r=i_0$, we 
obtain that
$$
\begin{array}{rcl}
{\displaystyle 
\left| \int_{0}^{2\pi } (r  \partial_r \varphi_n(q_{r}(\theta)) d\theta \, \right| } 
 & \leq & \sqrt{2}\mu^{\frac{1}{2}}_{g_n}(D_{r}(p))(64\mu_{g_n}^2(D_r(p)) + 
\mu_{g_n}(D_r(p))e_{g_n}(D_{r}(p)))^{\frac{1}{2}} \\ 
& \leq  & \sqrt{2}\mu^{\frac{1}{2}}_{g_n}(D_{r}(p))(64C^2_1+C_1C_2)^{\frac{1}{2}}\, .
\end{array}
$$
Thus, if 
$$
\psi_{n}(r)=\frac{1}{2\pi} \int_0^{2\pi} \varphi_n(r\cos{\theta},r\sin{\theta}) 
d\theta \, ,
$$
we have that
$$
\mid \psi_n(r)-\psi_n(0)\mid \leq \frac{1}{2\pi} \int_0^r \left| \int_0^{2\pi}
u \partial_u \varphi_n  d\theta \right| \frac{du}{u} \leq    
\frac{1}{2\pi} \int_0^r 
(2(64C^2_1+C_1C_2))^{\frac{1}{2}} \frac{\mu^{\frac{1}{2}}_{g_n}(D_{u}(p))}{u}du 
\, .
$$
 
By Lemma \ref{l55}, there exists a subsequence of the $\psi_n$s that has bounded
limit as $n\rightarrow \infty$. We relabel this subsequence as $\psi_n$ itself. 
By Lemma \ref{l56}, the singularity of $\mu_{g_n}(D_u(p))/u$ at $u=0$ is 
integrable, and the integral is uniformly bounded. 
These facts yield a contradiction since 
$\psi_n(0)=\varphi_n(0)\rightarrow \infty$.
\qed

\begin{lemma}
Let $\{g_n\}\in \mc{M}^{C_1,C_2}$ be an energy minimizing sequence in
the conformal class of the metric $g$ on $M$.
Suppose that $C_3$ is a constant such that
$$
\sup_{n} \sup_{q\in D(p)} \varphi_n (q) < C_3 \, . 
$$
Then for any relatively compact subdomain $\Omega \subset D(p)$, there exist 
a constant $C$ and a constant $\beta=\beta(C_1,C_2,C_3)$ such that
$$
\sup_{\Omega}\varphi_n  < \beta \inf_{\Omega} \varphi_n +C \, . 
$$
\end{lemma}

{\it Proof}. By the upper bound on the supremum of $\varphi_n$ and the $L^2$-bound
on the scalar curvature, it follows that $\varphi_n \in L^2(D(p))$ relative to the
flat metric, with a uniform bound on the $H^2$-norm.
If we write $\varphi_n= u_n+v_n$ with $\Delta_0 \varphi_n=\Delta_0 u_n$
in the interior of $D(p)$, and $\varphi_n \mid_{\partial D(p)}=v_{n}\mid_{\partial 
D(p)}$, by elliptic regularity for the Laplacian we see that $u_n \in H^2(D(p))$ 
uniformly, and by the Sobolev embedding theorem, its modulo is bounded by a constant
$C$ independent of $n$. It follows that the harmonic function
$v_n=\varphi_n-u_n$ is uniformly bounded above, and so by Schauder estimates, given 
any relatively compact subdomain $\Omega \subset D(p)$, there exist 
a constant $C$ and a constant $\beta=\beta(C_1,C_2,C_3)$ in $(0,1)$ such that
$$
\sup_{\Omega} (C-v_n) \leq \frac{1}{\beta}\inf_{\Omega} (C-v_n) \, .
$$
The desired result follows.
\qed
      
\begin{proposition}
Let $\{g_n=e^{\varphi_n}g\}\in \mc{M}^{C_1,C_2}$ be an energy minimizing 
sequence in the conformal class of the metric $g$ on $M$. Then for any 
subsequence of $g_n$, the embedded disk 
$D(p)$ carries at most $(64C^2_1+C_1 C_2)^{\frac{1}{2}}/(2\sqrt{2}\pi)$ 
bubble points, and there exists a subsequence that has
at most finitely many bubble points and no additional basic bubble points
in $D(p)$.
\end{proposition}

{\it Proof}. By Lemma \ref{le54}, a basic bubble point requires the concentration 
of a definite amount of area at the point, and so, there can be only finitely 
many of them. Suppose that $p_1, \ldots, p_k$ is an enumeration of all of these 
points. We have that 
$$
\sum_{j=1}^k A_{p_j} \leq \int_{D(p)} e^{\varphi_n}d\mu_0 \leq C_1\, ,
$$
and, by the estimate in Lemma \ref{le54}, that
$$
\frac{8\pi^2}{A_{p_j}} \leq 64A_{p_j}+E_{p_j} \, .
$$
Then,
$$
\begin{array}{rcl}
2(64C_1+C_2) & \geq & 2\sum_{j=1}^k (64 A_{p_j}+E_{p_j}) \\
  & \geq & {\displaystyle \sum_{j=1}^k \frac{16\pi^2}{A_{p_j}} } \vspace{1mm} \\ 
  & \geq & {\displaystyle \frac{16k^2 \pi^2}{\sum_{j=1}^k A_{p_j}}}\vspace{1mm}  \\   
  & \geq & {\displaystyle \frac{16k^2 \pi^2}{C_1}} \, ,   
\end{array}
$$ 
and the desired upper bound for $k$ follows.

If the sequence $g_n$ has $l$ distinct bubble points, and an additional basic bubble
$p$, by passing to a subsequence we make of $p$ a bubble point, and the subsequence
has now an additional bubble. This process must terminate in finitely
 many steps, or
else we would exceed the bound above for $k$.   
\qed

We let $\overline{H}^2(D(p))=H^2(D(p))\cup \{ \varphi_{-\infty} \}$, where 
$H^2(D(p))$ is the Sobolev space of order two on the domain $D(p)$ relative to
the flat metric, and $\varphi_{-\infty}=-\infty$ is the function that gives 
rise to the trivial tensor $g_{-\infty}=e^{\varphi_{-\infty}}g=0$. 
We can now state the two key consequences of these results in our context.
Their proof can be carried out exactly as in the proofs of 
\cite[Theorems 1, 3]{xxc}, using the lemmas above instead of their analogues
in \cite{xxc}.

\begin{theorem}
Let $\{g_n\}\in \mc{M}^{C_1,C_2}$ be an energy minimizing sequence in
the conformal class of the metric $g$ on $M$. Let $D(p)$ be an embedded 
disk in $M$ centered at $p$ as above where the metric $g_n$s can be 
represented as $g_n=e^{\varphi_n}g_0$, $g_0=dx^2+dy^2$ the standard flat 
metric in conformal coordinates on $D(p)$. Then there exists a subsequence
$\{ g_{n_j}\}$ with at most finitely many bubble points $p_1, \ldots, p_k$ in 
$D(p)$, $1\leq k \leq (64C^2_1+C_1 C_2)^{\frac{1}{2}}/(2\sqrt{2}\pi)$, and a 
metric $\tilde{g}=e^{\tilde{\varphi}}g_0$, 
$\tilde{\varphi}\in \overline{H}^{2}_{loc}(D(p)\setminus \{ p_1, \ldots, p_k\})$,
such that
$$
\varphi_{n_j} \rightarrow \tilde{\varphi} \; {\rm in} \; 
\overline{H}^{2}_{loc}(D(p)\setminus \{ p_1, \ldots, p_k\})\, ,
$$
and if the area and energy concentrations at $p_i$ are $A_{p_i}$ and
$E_{p_i}$, respectively, we have that  
$$
\lim_j \mu_{\varphi_{n_j}}(D(p))=\mu_{\tilde{\varphi}}(D(p)
\setminus \{p_1, \ldots, p_k\}) +\sum_{i=1}^k A_{p_i} \, ,
$$
and
$$
\lim_j e_{\varphi_{n_j}}(D(p))=e_{\tilde{\varphi}}(D(p)\setminus 
\{p_1, \ldots, p_k\}) +\sum_{i=1}^k E_{p_i} \, ,
$$
respectively. We obtain a subsequence of metrics $\{ g_{n_j}\}$ with 
a finite number of bubble points $\{ q_1, \ldots, q_k\}\subset M$, 
$k\leq (64C_1^2+C_1C_2)^{\frac{1}{2}}/2\sqrt{2}\pi$, and a metric
$\tilde{g}=e^{\tilde{\varphi}}g$ such that
$g_{n_j}\rightarrow \tilde{g}$ in $\overline{H}^2_{loc}(M\setminus
\{q_1, \ldots, q_k\})$, and we have the relations
$$
\lim_j \mu_{\varphi_{n_j}}(M)=\mu_{\tilde{\varphi}}(M
\setminus \{q_1, \ldots, q_k\}) +\sum_{i=1}^k A_{q_i} \, ,
$$
and 
$$
\lim_j \mc{E}_{g_{n_j}}(M)=\mc{E}_{\tilde{\varphi}}(M\setminus 
\{q_1, \ldots, q_k\}) +\sum_{i=1}^k E_{q_i} \, .
$$
\end{theorem}

\begin{theorem} \label{th511}
Let $\{g_n\}\in \mc{M}^{C_1,C_2}$ be a sequence of energy minimizing metrics on 
the disk $D(p)$ centered at $p$, where the metric $g_n$s can be 
represented as $g_n=e^{\varphi_n}g_0$, $g_0=dx^2+dy^2$ the standard flat 
metric in conformal coordinates $z=x+iy$ on $D(p)$. 
Suppose that $p$ is the only bubble point in $D(p)$ for the $g_n$s, that
the area and energy concentrations at $p$ are $A_p$ and $E_p$, respectively,
and that there exists $\tilde{\varphi}$ such that $\varphi_n \rightarrow 
\bar{\varphi}$ in $\overline{H}^2_{loc}(D\setminus \{ p\})$. Then we can choose
a sequence $\{\varepsilon_n \searrow 0\}$ so that the renormalized conformal 
factors $\phi_n=\varphi_n (\varepsilon_n(x,y))+\ln{\varepsilon_n}$ admit a 
subsequence $\{ \phi_{n_j}\}$ with finitely many bubble points $q_1, \ldots, q_k$,
$1\leq k \leq (64C^2_1+C_1 C_2)^{\frac{1}{2}}/(2\sqrt{2}\pi)$, and a conformal
factor $\tilde{\phi} \in \overline{H}(\mb{S}^2 \setminus 
\{ \infty, q_1, \ldots, q_k\})$ such that 
$$
\phi_{n_j}\rightarrow \tilde{\phi} \; \text{in $\overline{H}_{loc}
(\mb{S}^2 \setminus \{ \infty, q_1, \ldots, q_k\})$}\, ,
$$
and we have the relations
$$
A_p \geq \mu_{\tilde{\phi}}(\mb{S}^2 \setminus \{q_1, \ldots, q_k\})
+\sum_{i=1}^k A_{q_i} \, ,
$$
and
$$
E_p\geq \mc{E}_{\tilde{\phi}}(\mb{S}^2\setminus \{q_1, \ldots, q_k\})
+\sum_{i=1}^k E_{q_i} \, ,
$$
respectively. The renormalize sequence of metrics are of the form 
$\tilde{g}_n(z)=g_n(\varepsilon_n z+z(p_n))$ where $p_n \rightarrow p$, 
$p_n\in D(p)$ a point where the supremum of the area of $g_n$ in $D(p)$ occurs.
\end{theorem}
 
We may now prove a thin-thick decomposition for the limit of any
energy minimizing sequence in the conformal class of $(M_n,g_n)$. This 
is the analogue of the Cheeger-Gromov theorem \cite{cg} proven under a weaker
integral condition on the curvature tensors. 

\begin{theorem} \label{th512}
Let $(M_n,g_n)$ be any element of a $\mc{T}$-minimizing sequence 
$\{ (M_n,g_n)\}$ satisfying conditions {\rm M1-M5}. We set $C_1=2\pi(d(d-1)+2)$,
let $C_2$ be any uniform upper bound of all the $\mc{E}_{g_n}(M_n)$s, and  
consider any locally convergent sequence of energy minimizing metrics 
$\{g^n_m\}_m$ in the conformal class of $g_n$. Given any $\varepsilon >0$, 
there exist integers $N_{thin}$ and
$N_{thick}$,  depending only on $\varepsilon$ and $(64 C_1^2 +C_1 C_2 
)^{\frac{1}{2}}(2\sqrt{2}\pi)$, 
and a decomposition of $(M_n,g_m^n)$ into thick and thin components
$$
M_n =\cup_{j=1}^{N_{thick}} M_j^{thick}\cup \cup_{j=1}^{N_{thin}}M_j^{thin}\, ,
$$
such that 
\begin{enumerate}
\item Each $(M_j^{thick}, g_m^n\mid_{M_j^{thick}})$ converges in
$H^2_{loc}$ to a metric on $\mb{S}^2$ with a finite number of small discs 
deleted whose sizes can be made as small as possible.
\item Each $M_j^{thick}$ is connected, but no any two of them forms a 
connected subset of $M_n \setminus \cup_{j=1}^{N_{thin}}M_j^{thin}$.
\item Each $(M_j^{thin},g_m^n\mid_{M_j^{thin}})$ is homeomorphic to a cylinder 
$\mb{S}^1 \times (a,b)$, and the length of any loop $\mb{S}^1 \times {p}$, 
$a<p<b$ is strictly less than $\varepsilon$.
\end{enumerate}
\end{theorem}  

\begin{corollary} \label{c513} 
Exactly $d$ of the thick components $M_j^{thick}$ are homeomorphic to 
a homologically nontrivial sphere in $\mb{P}^2(\mb{C})$ with a finite number of 
points deleted, and all of the remaining thick components as well as all of 
the thin components bound in $\mb{P}^2(\mb{C})$.
\end{corollary}

{\it Proof}. We can shrink each of the thin components 
to produce in the limit a homology representative of
$[S_d]$ in $\mb{P}^2(\mb{C})$
consisting of $2$-spheres any two of which intersect at a point.
For dimensional reasons, these spheres can be made to intersect transversally.
The result follows.
\qed

In the renormalization process of Theorem \ref{th511}, a neck is created 
for each of the bubbles of the initial one. The length of a circle in these
necks is bounded above by $\varepsilon$, and as $n$ increases to $\infty$, the
conformal distance between the bounding circles approaches $\infty$ while 
part of the interior of the neck collapses into a line. The collapsing occurs 
while either the scalar curvature or the diameter of the
neck remain bounded, and the size of the neck can be made arbitrarily small
by letting the lengths of the bounding circles go to zero.       
 
\begin{corollary} \label{c514}
There exists an energy minimizing sequence $\{g^n_m \}_m$ in the conformal 
class of $g_n$ such that $(M_n,g_m^n)$ converges 
to a transversal intersection of $d$ complex lines each with its induced 
Fubini-Study metric, and a collection of $2$-spheres intersecting transversally
each of which bounds in $\mb{P}^2(\mb{C})$, where the energy of the limit 
metric over the portion of the manifold that bounds is equal to the barrier 
$b_n$ of the metric $g_n$.
\end{corollary}

{\it Proof}. We deform conformally each of the $d$ nonbouding spheres into
complex lines, pushing their initial excess area through the necks onto the 
bounding spheres.  Complex lines are totally geodesic submanifolds of $\mb{P}^2(
\mb{C})$, and so they contribute nothing to the energy, which must then come
from the remaining portion of the limiting manifold. But the limiting manifold
realizes the infimum of the energy, which is given by the barrier $b_n$ of the
conformal class of $(M_n,g_n)$. 
\qed

\subsection{Precompactness properties of $\mc{T}$ over $\mc{M}_{[S_d]}(\mb{P}^2
(\mb{C}))$}

We now study the asymptotics of the functional $\mc{T}$ over representatives of 
$[S_d]\in H_2(\mb{P}^2(\mb{C}); \mb{Z})$ as these bubble off into the transverse 
intersections of $d$ complex lines that absolutely minimize the energy functional
over all conformal classes. For that, first we look in further detail at the 
resolution of singularities of the latter introduced in \S \ref{s51}. 

We consider an immersed complex singular representative $\Sigma $ of $[S_d]$
whose singularities are modeled in local complex coordinates $(z_1, z_2)$ by 
$$
z_1z_2=0 \, .
$$
If $z_1=w_1+iw_2$ and $z_2=w_1-iw_2$, then $z_1z_2=w_1^2+w_2^2$, 
and so the equation above is just the zero locus of the quadratic complex 
function $q(w)=w_1^2+w_2^2$. We regularize the singularity by considering 
$$
z_1z_2 =q(w)=w_1^2+w_2^2=\varepsilon 
$$
for $0<\varepsilon \ll 1$, and apply the surgery procedure discussed in 
\S \ref{s51} to the linked discs $L=\{ z=(z^1,z^2): \; z^1 z^2=0,\; |z|<r\}$ 
to obtain a representative 
$\Sigma'$ of the class $[\Sigma]$ that is also smooth in a neighborhood of $z=0$. 

Notice that if $w_j = x_j +i y_j$, $j=1,2$, $x_j,y_j$ real, then
the singular fiber $q(w)=0$ intersects the Euclidean sphere 
$\mb{S}_{r}^3= \{ w=(w_1, w_2) : \; |w|=r \}$ in the $r/\sqrt{2}$-Euclidean sphere bundle 
of the cotangent space $T^* \mb{S}_{r/\sqrt{2}}^1$ of the Euclidean circle
$\mb{S}_{r/\sqrt{2}}$,  for we must have the relations
$$
\begin{array}{rcl}
\sum_{j=1}^2 x_j^2-y_j^2 & =  & 0 \, , \\
\sum_{j=1}^2 x_jy_j & =  & 0 \, , \\
\sum_{j=1}^2 x_j^2+y_j^2 & =  & r^2 \, .
\end{array}
$$
Thus, the two loops $L=\{ z_1 z_2=0 \} \cap \mb{S}^3_r$ are linked, and have
linking number $1$, and the desingularization procedure replaces the intersecting
discs $\{ z_1 z_2=0\, : \; |z|\leq r\}$ by the minimal genus Seifert surface 
for the Hopf link, an annulus with boundary $L$. 

Suppose now that $\Sigma$ has singularities of the type above relative to 
an affine system $(z_1, z_2)$ of local coordinates in $\mb{P}^2(\mb{C})$. 
We look at the jump in the extrinsic scalar curvature of $\Sigma$ arising from
regularizing the discs $D_r=B_r(0)\cap \{ z_1z_2=0\}\subset \Sigma$ to 
$D_{r,\varepsilon}=B_r(0)\cap \{ z_1z_2=\varepsilon\}$ according to the 
procedure above. We denote by $\Sigma_{\varepsilon}$ the embedded representative
of $[\Sigma]$ so obtained.  
We study the quantities
$$
\Theta_r =\int_{D_r} 2K(e_1,e_2)d\mu_g = 8 \int_{D_r} d\mu_g \mid_{D_r}
$$
and
$$
\Theta_{r,\varepsilon} =\int_{D_{r,\varepsilon}} 2K(e_1,e_2)d\mu_g =
8 \int_{D_{r,\varepsilon}} d\mu_g \mid_{D_{r,\varepsilon}} \, , 
$$ 
respectively. We denote by 
$\Theta({\Sigma})$ and $\Theta(\Sigma_{\varepsilon})$ the
total extrinsic scalar curvatures of the minimal current $\Sigma$ and the
embedded submanifold $\Sigma_{\varepsilon}$,
respectively.

The Fubini-Study metric restricted to either one of the 
discs in $D_r$ is given by 
$$ 
\omega= \frac{i}{2(1+|z|^2)^2} dz \wedge d\bar{z} \, ,
$$
while its restriction to $D_{r,\varepsilon}$ is given by 
$$
\omega_{\varepsilon}=\frac{i}{2(1+|z|^2+\frac{\varepsilon^2}{|z|^2})^2}\left( 
1+4 \frac{\varepsilon^2}{|z|^2}+ \frac{\varepsilon^2}{|z|^4}\right) dz 
\wedge d\bar{z} \, .
$$
We then have that 
$$
\Theta_r =8 \int_{D_r} d\mu_g = 8\pi \frac{2r^2}{2+r^2}\, ,
$$
and  
$$
\Theta_{r,\varepsilon} =
8\int_{D_{r,\varepsilon}} \frac{i}{2(1+|z|^2+\frac{\varepsilon^2}{|z|^2})^2}
\left(1+4 \frac{\varepsilon^2}{|z|^2}+ \frac{\varepsilon^2}{|z|^4}\right) dz 
\wedge d\bar{z}= 8\pi + 8\pi \frac{r^2}{1+r^2} +o(\varepsilon) \, ,
$$
respectively.

\begin{lemma} \label{l515}
For the embedded representative $\Sigma_{\varepsilon}$ of $[\Sigma]$, we have
that 
$$
\lim_{\varepsilon}\mc{T}(\Sigma_{\varepsilon})=4\pi d(d-1)\, , \quad
\lim_{\varepsilon}\Psi(\Sigma_{\varepsilon})=0\, , \quad
\lim_{\varepsilon}\Theta(\Sigma_{\varepsilon})=
\Theta(\Sigma)=8\pi d\, .
$$
\end{lemma}

{\it Proof}. The minimal current $\Sigma$ is the weak limit of the embedded 
current $\Sigma_{\varepsilon}$, and so $\Psi(\Sigma_{\varepsilon})\rightarrow
0$ from above as $\varepsilon \rightarrow 0$. 
Letting $r\rightarrow \infty$ in the 
asymptotic expansions found above, we see that $\Theta(\Sigma_{\varepsilon})
-\Theta(\Sigma)$ approaches zero as $\varepsilon \rightarrow 0$.

The statement about $\mc{T}(\Sigma_{\varepsilon})$ then follows from the
identity
$$
\mc{T}(M)=\Theta(M)+\Psi(M)-\int_M s_g d\mu_g \, .
$$
We just need to observe that the regularization $\Sigma_{\varepsilon}$ of 
$\Sigma$ at each of its $d$-intersection points results into a surface of genus 
$\bn{n-1}{2}=\frac{(d-1)(d-2)}{2}$, and so applying the Gauss-Bonet theorem, 
we obtain that
$$
\mc{T}(M_{\varepsilon})=\Theta(M_{\varepsilon})+\Psi(M_{\varepsilon}
)-8\pi \left( 1-\frac{(d-1)(d-2)}{2}\right)\rightarrow 4\pi d(d-1) \, .
$$
This finishes the proof.
\qed  

We now have the following:
\begin{theorem} \label{th516}
Over $\mc{M}_{[S_d]}( \mb{P}^2(\mb{C}))$,
the functionals {\rm (\ref{cpp})} is globally 
bounded below,
$$
{\mathcal M}_{[S_d]}(\mb{P}^2(\mb{C}))\ni \Sigma
\mapsto \mc{T}(\Sigma)=\int_{\Sigma} \tg(\a,\a) d\mu_g 
\geq 4\pi d(d-1) \, ,
$$
and the lower bound is achieved if, and only if, $\Sigma=S_d$ is an algebraic 
curve of degree $d$.
\end{theorem}

{\it Proof}. Consider a minimizing sequence $(M_n, g_n)$ satisfying M1-M5. 
By M3, the topology of all the $M_n$s is fixed and equal to the topology of 
a fixed surface $\Sigma$, and by Theorem \ref{th512}, we get a decomposition of
$(M_n,g_n)=(f_n(\Sigma),g_n)$ as a union of a fixed number of 
topological spheres tubed together by cylinders, with exactly $d$ of the 
spheres nontrivial in homology. It follows that $t$ cannot strictly less than
$4\pi d(d-1)$. For it if were, then then from some sufficiently large $N$ on, 
the barrier $b_n=\mc{T}^2(M_n)/\mu_{g_n}$ would be strictly less than the barrier
of a complex curve of degree $d$, which cannot be by 
Corollary \ref{c514} and Lemma \ref{l515}.
Since each of the bounding spheres can be 
collapsed into a cylinder with negligible contribution to the energy and to
$\mc{T}$ as $n\rightarrow \infty$, at $\mc{T}(M_n)\rightarrow t$, 
the areas $\mu_{g_n}$ must converge to $\pi d$, the area of a complex curve of
degree $d$ in $\mb{P}^2(\mb{C})$. It follows then that if $(M_n,g_n)$ develops
singularities as $n\rightarrow \infty$, the value of $\mc{T}(M_n)$ becomes
infinitesimally close to $4\pi d(d-1)$ but not quite exactly that value. On
the other hand, by Proposition \ref{p5}, the value of $\mc{T}$ for any
complex curve of degree $d$ is exactly $4\pi d(d-1)$. Thus, the  
infimum $t$ is achieved, and the minimizing sequence $(M_n,g_n)$ either 
converges to a curve of degree $d$ in $\mb{P}^2(\mb{C})$, or it degenerates
into an intersection of $d$-complex lines union perhaps additional line segments
tying two of these lines together, segments that 
could be homotoped away to a point.  
\qed

\begin{theorem}
{\rm (Kronheimer-Mrowka \cite{krmr})} Let $\Sigma$ be a connected 
oriented smooth 
embbeded surface in ${\mathbb C}{\mathbb P}^2$ that represents the 
homology class of an algebraic curve of degree $d$. If $g_{\Sigma}$ is its
genus, we have that
$$
g_{\Sigma} \geq \frac{(d-1)(d-2)}{2}\, ,
$$
and the equality is achieved if $\Sigma $ is a complex
curve in ${\mathbb C}{\mathbb P}^2$.
\end{theorem}

{\it Proof}. If there were an embedded representative of $S_d$ of genus less
than $\bn{n-1}{2}$, then by minimizing the energy in the conformal class of such
a representative, and then collapsing the spheres that bound produced by
Theorem \ref{th512}, we would obtain a representative of the class $[S_d]$ with
a barrier strictly less than the barrier of a complex curve of degree $d$, 
which cannot be.
\qed 
  

\section{Homotopic submanifolds of varying critical levels for $\mc{S}$}
\label{s6}
\setcounter{theorem}{0}
In their alternative to compactification, Randall and Sundrum \cite{rasu}
consider a three brane in a five dimensional bulk space, and quantize the 
system and treat the non-normalizable modes by introducing a second brane 
at a certain distance from the first. Brane directions and time yield
space-time slices. The combined branes are to be placed at the boundaries
of a fifth dimension. The action function they consider is the sum of 
suitable multiples of the volume and total scalar curvature in space-time,
together with boundary contributions from the branes. Then they find a 
warped product solution to Einstein's equation on the five dimensional 
bulk that carries two space-time slices that are invariant under the 
Poincar\'e group, the effective action at each of these being different from
each other. In this section we use the functional (\ref{ssf}) to exhibit 
examples with the same flavour as that exhibited by the work 
of Randall and Sundrum.

Indeed, we study critical points of (\ref{ssf}) that are 
{\it localized} in a particular direction, and such that 
two different critical points can be homotopically
deformed into each other through embedded submanifolds. 
In our examples, the fixed backgrounds are Riemannian manifolds
of dimension four, but we can take 
their metric product with a fifth time direction. After doing so, the
said critical points multiplied by time 
become the analogues of the brane-time slices 
of Randall and 
Sundrum, and the different critical values correspond to these slices'
different effective actions.
In each of our last two examples here, the second critical point occurs 
asymptotically at $\infty$ and is flat. The ``nonasymptotic'' 
critical points of all of our examples are neither
critical points of (\ref{sf}) nor of (\ref{mc}) separately.

By Theorem \ref{th2}, the critical point equation of
(\ref{ssf}) in the codimension one case is given by
\begin{equation}
2 k_i K_{\tg}(e_i,\nu)-2h \sum_i K_{\tg}(e_i,\nu)-2e_i(r_{\tg}(e_i,\nu))
+2{\rm trace}A_\nu^3 -
3\|\alpha\|^2 h + h^3 = 0 \, ,
\label{c1c}
\end{equation}
where $\nu$ is a normal field along $M\subset \tm$, $\{ e_1, \ldots, e_n\}$
is the orthonormal frame of principal directions with associated principal
curvatures $k_1, \ldots, k_n$, and $r_{\tg}$ is the Ricci curvature tensor
of $\tg$. If the background manifold is 
Einstein, this critical point equation simplifies to 
\begin{equation}
2\sum_i  k_i K_{\tg}(e_i,\nu)-2h \sum_i K_{\tg}(e_i,\nu)+2{\rm trace}A_\nu^3 -
3\|\alpha\|^2 h + h^3 = 0 \, .
\label{c2c}
\end{equation}

In our examples with two isotopic solutions to this equation,
the {\it localization} occurs in the normal direction.

\subsection{Cartan's isoparametric families}
Let us reconsider the isoparametric families $M_t^{2n}\subset 
{\mathbb S}^{2n+1}$ of Example \ref{ex8} (re)discovered
by Nomizu \cite{no}. The manifold $M_t^{2n}$ has principal curvatures 
$k_1(t)$ and $k_2(t)$ each with multiplicity one, and $k_3(t)$ and $k_4(t)$ 
each with multiplicity $(n-1)$. Thus, the mean curvature function
is given by $h(t)= k_1(t)+k_2(t)+(n-1)(k_3(t)+k_4(t))$, 
while $a(t)={\rm trace}A_\nu^3=k^3_1(t)+k^3_2(t)+(n-1)(k^3_3(t)+k^3_4(t))$,
and $M_t^{2n}$ is a critical point of
(\ref{ssf}) if, and only if, $t\in (0,\pi/4)$ is a root of the equation
$$
2(1-2n)h(t)+2a(t)
-3h(t)((k^2_1(t)+k^2_2(t)+(n-1)(k^2_3(t)+k^2_4(t)))+
h^3(t)=0 \, .
$$

For $n=1$ this equation has no solution on $(0,\pi/4)$. For $n=2$, there
is only one solution in the indicated range, and the corresponding 
critical submanifold $M_t$ is the austere hypersurface 
of ${\mathbb S}^5$ already pointed out
in Example \ref{ex8}. The critical value is $12$ times the
volume of $M_t$. For $n=3$, there is also exactly one solution 
in the said range, $t=0.5268183350$, but the corresponding critical point 
$M_t^6$ is not
a minimal hypersurface of ${\mathbb S}^{7}$. The critical value is 
$19.71086118$ times the volume of $M_t$.

For $n\geq 4$ the equation above has two solutions in the range
$(0,\pi/4)$, and their critical values are real numbers of opposite signs.
For instance, for $n=4$ (so the brackground bulk is of dimension 10, space and
time together), one root is $t=0.1830436696$, and the 
corresponding critical value for the 
hypersurface $M_t^8\subset {\mathbb S}^9$ is $-131.2969104$ times the volume.
The other root occurs at $t=0.5770248421$, and the critical value of the
associated $M_t^8$ is $27.29691039$ times the volume. Each of  
these two critical points $M_t^{8}$ inherits 
a metric of positive scalar curvature.

\begin{remark}
The example above is unsatisfactory in that the homology classes 
represented by the exhibited critical submanifolds are trivial. 
Examples otherwise 
would require a background with nontrivial homology.
\qed
\end{remark}

\subsection{The Schwarzschild background}
The Schwarzschild metric is a Ricci flat metric on the total space of the
trivial
bundle ${\mathbb S}^2 \times {\mathbb R}^2 \rightarrow {\mathbb S}^2$. It
can be cast as a doubly warped product on $(-\infty,\infty)\times 
{\mathbb S}^1\times {\mathbb S}^2$, 
$$
\tg= dr^2 + \varphi^2(r) ds_{1}^2 + \psi^2(r) ds_2^2 \, ,
$$
with suitable asymptotics at $r=0$ and $r=\infty$, respectively.    
We have that
\begin{equation}
\dot{\psi}^2 = 1-\frac{\beta}{\psi} \, , \; \; 
\dot{\psi}=\frac{1}{2\beta} \varphi \, , \label{ode}
\end{equation}
where $\psi(0)=\beta >0$. The solution to the equation satisfied by 
$\psi$ exists for all $r\in (-\infty,\infty)$; it is an even function
of $r$, 
and $\psi(r)/|r| \rightarrow 1$ as $r\rightarrow \pm \infty$.

We take $\nu= \partial_r$, the gradient of the distance function $r$, 
and $\{ e_1, e_2, e_3\}$ an orthonormal set
with $e_1$ tangent to the ${\mathbb S}^1$ factor, and $e_2$ and $e_3$ tangent
to the ${\mathbb S}^2$ factor.  The principal curvatures of the
${\mathbb S}^1 \times {\mathbb S}^2$ hypersurfaces defined by 
setting $r$ to a constant are 
$$
k_1 = -\frac{\dot{\varphi}}{\varphi}\, , \; \; k_2 =k_3=-
 \frac{\dot{\psi}}{\psi}\, .
$$
We have that
$$
K_{\tg}(e_1, \nu) = -\frac{\ddot{\varphi}}{\varphi}\, ,  \; \;
K_{\tg}(e_2, \nu) = K_{\tg}(e_3, \nu) = -\frac{\ddot{\psi}}{\psi} \, .
$$
Notice that $\sum_i K_{\tg}(e_i, \nu)=0$ 
follows on the basis of these expressions for the
sectional curvatures, and (\ref{ode}), a fact consistent with the Ricci
flatness of the metric.

The critical hypersurfaces slices $r=r_0$ for (\ref{ssf}) are determined by 
the roots of the function 
$$
%
C= 2\left(\frac{\dot{\varphi}}{\varphi}\frac{\ddot{\varphi}}{\varphi}
+2\frac{\dot{\psi}}{\psi}\frac{\ddot{\psi}}{\psi}\right)  -6 
\frac{\dot{\varphi}}{\varphi}\left(\frac{\dot{\psi}}{\psi}\right)^2 \, .
$$
This is an odd function of $r$ that is positive to the immediate
right of $0$. Using the asymptotic behaviour of $\psi$, we see that
it becomes negative at some point, and so it must have a
a zero $r_0>0$, which yields a critical slice ${\mathbb S}^1\times
{\mathbb S}^2$ for the functional (\ref{ssf})
with metric $\varphi^2(r_0)ds_1^2+\phi^2(r_0)ds_2^2$, the 
${\mathbb S}^1$-factor being flat while the ${\mathbb S}^2$-factor is
Einstein of constant $(-\ddot{\psi}/\psi+(1-\dot{\psi}^2)/\psi^2-
\dot{\varphi}\dot{\psi}/\varphi\psi)\mid_{r=r_0}$. Thus, this 
critical slice is provided with a geometric structure modeled by
$\mb{R}\times \mb{S}^2$, one of the eight such structures
in Thurston's geometrization program. 

Further, as $r\rightarrow \infty$, the function $C$ increases up to $0$, and
we obtain asymptotically a second critical point for (\ref{ssf}) that
looks like ${\mathbb S}^1 \times {\mathbb R}^2$ with the flat product 
metric, and where the metric in the first factor is $(2\beta)^2 ds_1^2$.
The geometric model of this critical point is $\mb{R}^3$. 

The critical hypersurface at 
small distance is highly curved in comparison with 
the critical point at large distance, which is a scalar flat space. 
This behaviour is reminescent of some of the features of
{\it spacetime foam}, which is taken to be the building block of the universe 
in some theories of gravity. Notice that 
our functional (\ref{ssf}) is just a linear combination of the total scalar 
curvature and a functional equivalent (in the sense that it is bounded above 
and below by a suitable multiple of it) to the volume form of the submanifold,
with the latter being  exactly a constant times the volume 
in the case of hypersurfaces defined by setting $r$ equal to a constant.

By the Ricci flatness of the background metric, (\ref{sc}) implies that
 the scalar curvature of the critical slice determined
by $r=r_0$ is given by the value of the function
$$
S=2\frac{\dot{\psi}}{\psi}\left( \frac{\dot{\psi}}{\psi}+2\frac{\dot{\varphi}}
{\varphi}\right)
$$
at $r=r_0$. Thus, the first of the critical points has positive scalar 
curvature. The asymptotic critical point at $\infty$ has zero scalar
curvature. In fact, $S(r)\approx (2/r)(1/r+\beta/r^2)$ as $r\rightarrow
\infty$.

\subsection{The Eguchi-Hanson background}
The Eguchi-Hanson metric is a family of scaled 
Ricci flat metric on the total space of the tangent
bundle $T{\mathbb S}^2\rightarrow {\mathbb S}^2$ to the $2$-sphere. If 
we use the coframe $\{ \sigma^1, \sigma^2, \sigma^3\}$ 
that identifies ${\mathbb S}^3$ as the Lie group ${\mathbb S}{\mathbb U}(2)$,
the metric can be written as a doubly warped product 
$$
\tg=dr^2 + \varphi^2(r)(\psi^2 (r)(\sigma^1)^2 +(\sigma^2)^2+(\sigma^3)^2)\, ,
$$
where 
$$
\dot{\varphi}=\psi\, , \; \;  \dot{\varphi}^2 = 1-k\varphi^{-4} \, ,
$$
subject to the boundary conditions $\varphi(0)=k^{\frac{1}{4}}$, 
$\dot{\varphi}(0)=0$, $\psi(0)=0$, and $\dot{\psi}(0)=2$. The parameter 
$k$ can be changed via scaling.

We take $\nu= \partial_r$, the gradient of the distance function $r$, 
and $\{ e_1, e_2, e_3\}$ dual to $\{\varphi \psi \sigma^1, 
\varphi \sigma^2, \varphi \sigma^3\}$.
The principal curvatures of the hypersurfaces defined by 
setting $r$ to a constant are given by
$$
k_1 = -\frac{\partial_r(\varphi \psi )}{\varphi \psi}\, , \; \; k_2 =k_3=-
 \frac{\partial_r{\varphi}}{\varphi}\, ,
$$
while
$$
K_{\tg}(e_1, \nu) = -\frac{\partial_r^2 (\varphi \psi)}{\varphi \psi}\, ,
 \; \;
K_{\tg}(e_2, \nu) = K_{\tg}(e_3, \nu) = -\frac{\partial_r^2{\varphi}}
{\varphi} \, .
$$

The critical hypersurfaces slices $r=r_0$ for (\ref{ssf}) are determined 
then by the roots of the function 
$$
C= 2\left( \frac{\partial_r(\varphi \psi)}{\varphi \psi}
\frac{\partial_r^2(\varphi \psi)}{\varphi \psi}
+2\frac{\partial_r \varphi}{\varphi}\frac{\partial_r^2{\varphi}}{\varphi}
\right)  -6 \frac{\partial_r(\varphi \psi)}{\varphi \psi}\left(
\frac{\partial_r{\varphi}}{\varphi}\right)^2 \, .
$$
This function is positive to the immediate
right of $0$, and by the asymptotic behaviour of $\psi$, it becomes 
negative at some point, and increases to $0$ as $r\rightarrow \infty$.
So the function must have a
a zero $r_0>0$, yielding a critical slice 
for the functional (\ref{ssf})
with metric $\varphi^2(r_0)(\psi^2(r_0)(\sigma^1)^2+(
(\sigma^2)^2+(\sigma^3)^2)$, a rescaled Berger sphere, and an asymptotic
critical slice at $r=\infty$.  The asymptotic critical point is an  
$\mb{P}^3(\mb{R})$. The cone on it as a base has the asymptotic flat metric
$r^2((\sigma^1)^2+(\sigma^2)^2+(\sigma^3)^2)$, and the metric on the asymptotic
critical slice is induced by it. As in the previous example,
the critical hypersurface at 
small distance is highly curved, while the critical point at large distance 
is a scalar flat space.

The scalar curvature of a critical slice in the metric induced by the background
Eguchi-Hanson metric is given by the value of the function
$$
S=2\frac{\partial_r{\varphi}}{\varphi}\left( \frac{\partial_r{\varphi}}{
\varphi}+2\frac{\partial_r(\varphi \psi)}
{\varphi \psi}\right)
$$
at $r=r_0$. Thus, the first of the critical points has positive scalar 
curvature, while the asymptotic critical point at $\infty$ has zero scalar
curvature, as already observed.
In fact, $S(r)\approx (2/r)(3/r+4k/r^5)$ as $r\rightarrow \infty$.

\subsection*{Acknowledgement}
We would like to express our gratitude to Professor Wojciech Kucharz for
several conversations that took place during during the preparation of this 
paper, and from which we benefitted substantially.

\end{document}